\documentclass[12pt,1eqno]{amsart}
   \usepackage{amsthm}
   \usepackage{amsmath}
   \usepackage{amssymb}
   \usepackage{amsfonts}
   \usepackage{amscd}
   \usepackage[mathscr]{eucal}
   \usepackage{color}
   \usepackage{verbatim}
   \usepackage{graphics}
  \usepackage{hyperref}

\newcommand{\bbb}{\mathcal{B}}
\newcommand{\ccc}{\mathcal{C}}

\newcommand{\fff}{\mathcal{F}}

\newcommand{\llL}{\mathcal{L}}

\newcommand{\wte}{\widetilde{E}}

\newcommand{\wtec}{\widetilde{E^c}}
\newcommand{\wtx}{\widetilde{x}}
\newcommand{\wt}{\widetilde}

\newcommand{\con}{\equiv}
\newcommand{\m}{\pmod{m}}

\newcommand{\zp}{\mathbb{Z}/p\mathbb{Z}}
\newcommand{\zq}{\mathbb{Z}/q\mathbb{Z}}
\newcommand{\zm}{\mathbb{Z}/m\mathbb{Z}}

\newcommand{\Card}{\operatorname{Card}}

\newcommand{\tp}{\tilde{\pi}_2}

\newtheorem{thm}{Theorem}[section]
\newtheorem*{mlem*}{Main Lemma}
\newtheorem*{thm*}{Theorem}
\newtheorem{prop}[thm]{Proposition}

\newcommand{\Z}{{\mathbb Z}}

\newtheorem{lem}[thm]{Lemma}

\newtheorem*{cor*}{Corollary}
\newtheorem*{conj*}{Conjecture} 
\newtheorem{cor}[thm]{Corollary}

\theoremstyle{definition}

\newtheorem{rem}[thm]{Remark}

\begin{document}

\title{Sum-free subsets of finite abelian groups of type III}
\author{R. Balasubramanian, Gyan Prakash and D.S. Ramana\\  %The Institute of Mathematical Sciences\\ 
%CIT Campus, Taramani\\ %Chennai-600113, India.\\ %E-mail:
%\texttt{\{balu,gyan\}@imsc.res.in}  
} 
%\address{Institute of Mathematical Sciences,
%CIT Campus, Tharamani, Chennai 600 113, India}
%\email[R.Balasubrmanian]{balu@imsc.res.in}
%\address{Harish-Chandra Research Institute, Chhatnag Road, Jhunsi, Allahbad-211019, India}
%\email[Gyan Prakash]{gyan@hri.res.in, gyan.jp@gmail.com}
%\email[D.S. Ramana]{suri@hri.res.in}
\date{} 
%\maketitle
%\tableofcontents

\begin{abstract}
\noindent
A finite
abelian group $G$ of cardinality $n$ is said to be of type III if every prime divisor of $n$ is congruent to 1 modulo
3. We obtain a classification theorem for sum-free subsets of 
largest possible cardinality in a finite abelian group $G$ of type III. This theorem, when taken together with known results, gives a complete
characterisation of sum-free subsets of the largest cardinality in any finite
abelian group $G$. We supplement this result with a theorem on the structure of sum-free subsets of  cardinality ``close" 
to the largest possible in a type III abelian group $G$. We then give two applications of these results. Our first
application allows us to write down a formula for the number of orbits under
the natural action of ${\rm Aut}(G)$ on the set of sum-free subsets of $G$ of
the largest cardinality when $G$ is of the form $({\Z}/m{\Z})^r$, with
all prime divisors of $m$ congruent to 1 modulo 3, thereby extending a result of
Rhemtulla and Street. Our second application provides an upper
bound for the number of sum-free subsets of $G$. For finite abelian groups $G$
of type III and with {\em a given exponent} this bound is substantially better
than that implied by the bound for the number of sum-free subsets in an
arbitrary finite abelian group, due to Green and Ruzsa. 
\end{abstract}
\maketitle
\section{Introduction}

\noindent
A subset $A$ of an abelian group $G$ is said to be {\em sum-free}
if the sum of any pair of elements of $A$ lies in the complement of $A$ { in $G$}. In
 other words, $A$ is sum-free if there is no solution to the equation $x+y
= z$ with $x, y, z$ in  $A$. For example, the set of odd integers
is a sum-free subset of the group ${\Z}$.

\vspace{2mm}
\noindent
%The study of sum-free subsets of abelian
 %groups is the subject of a number of investigations in the
 %literature, the most relevant to this article being the
 %work~\cite{GR} of Ben Green and Imre Ruzsa. 
 
 \vspace{2mm}
 \noindent
 For the purpose of
 studying their sum-free subsets, it is convenient to classify finite
 abelian groups into three {\em types} (see Subsection \ref{2.5} below). A finite abelian group $G$  is said to be of type III if every prime divisor of ${\rm Card}(G)$ is congruent to $1$ mod  $3$. 
 
 %The present article is concerned with  
 %sum-free subsets in a finite abelian group of type~III.  

 %In particular, we obtain a structure of
 %sum-free subsets (Theorem~\ref{clas}) of the largest possible cardinality and also give a count for the total number of
 %sum-free sets (Theorems~\ref{orbi} and~\ref{coun}) in a finite
 %abelian group of type III.

\vspace{2mm}
\noindent
 %When $G$ is not of type III, 
 Let us write $c(G)$ for the largest possible cardinality of the 
 of a sum-free set in a given abelian group $G$. It is known that if $G$ is of exponent $m$ and cardinality $n$ then

%\begin{thm}\cite[Theorem 1.5]{GR} \label{mug}
%$When $G$ is a finite abelian group of exponent $m$ and cardinality $n$,
%the largest possible cardinality $c(G)$ of a sum-free subset in $G$ is equal to
\begin{equation}\label{mug}
c(G) = n\left(\max_{d | m} \frac{\lfloor\frac{d-2}{3}\rfloor +1}{d}\right).
\end{equation}
%\end{thm}

\noindent
This was first shown for abelian groups $G$ that are not of type III by P.H. Diananda and
 H.P. Yap~\cite{DY}. The relation \eqref{mug} was confirmed for abelian groups $G$ of type III only relatively recently  by  B.J. Green and I. Ruzsa in their important work \cite{GR}.

\vspace{2mm}
\noindent
The main result of the first part of this article is a classification theorem,  Theorem~\ref{clas}, which characterises sum-free subsets of cardinality $c(G)$ in finite abelian groups $G$ of type III. When this theorem is read together with the results of Diananda and Yap~\cite{DY}, one obtains a complete characterisation of sum-free subsets of the largest possible cardinality in all (three) types of finite abelian groups.

\vspace{2mm}
\noindent
We state Theorem \ref{clas} with the aid of the following notation. We recall that when $G$ is an abelian group and $M$ is a subgroup of $G$, a {\em supplement}
of $M$ in $G$ is a subgroup $H$ of $G$ such that the
canonical map $(x,y) \rightarrow x+y$ from $H \oplus M$ into $G$ is an
isomorphism. Let us now suppose that $G$ is a finite abelian group of exponent $m$. Then $G$ contains a subgroup isomorphic to $\zm$ and every
such subgroup of $G$ has a supplement in $G$. By a {\em splitting} of $G$ by ${\Z}/m{\Z}$ we shall mean a pair $(H,f)$
where $f$ is an injective homomorphism from $\zm$ into $G$ and $H$ is
a supplement of the image of $f$ in $G$. Given a splitting $(H,f)$ 
of $G$ by ${\Z}/m{\Z}$ and $B$, $C$ subsets
of $H$ and ${\Z}/m{\Z}$ respectively, we shall write $(B,C)_{(H,f)}$, or simply
$(B,C)$, to denote the subset $B + f(C)$ of $G$. Moreover, for any 
$x_1 \in H$ and  $x_2 \in \Z/m\Z$ we write $(x_1,x_2)$ to denote the
element $x_1+f(x_2)$ of $G$. It is easily seen that given
for each $x\in G$, there { exists} a unique $x_1 \in H$ and {a} unique
$x_2 \in \Z/m\Z$ such that $x= (x_1,x_2).$

\vspace{2mm}
\noindent
When $X$ is a subset of ${\Z}$ and $m$ is an integer we write $X_m$ to denote the canonical image of $X$ in ${\Z}/m{\Z}$. For example, suppose that $a$ and $b$ are integers with $a < b$. We then write $[a,b]_m$ to denote the canonical image in $ \Z/m\Z$ of the set of integers in the real interval $[a,b]$. 

\begin{comment}
\vspace{2mm}
\noindent
Finally, let us note that since every divisor of the cardinality of a finite abelian group $G$ of type III is congruent to 1 modulo 3, the exponent $m$ of $G$ is odd and therefore it is congruent to 1 modulo 6.  
\end{comment}

\vspace{3mm}
\begin{thm}\label{clas}
Suppose that $G$ is a finite abelian group of type III and exponent $m$. Let $l$ denote $\frac{m-1}{3}$ and let $(H,f)$ be a splitting of $G$ by $\zm$. Also, let $K$ be a subgroup of $H$ and $K^c$ denote the complement of $K$ in $H$. Then each of the following is a sum-free subset of $G$ of the largest possible cardinality. 
 
\vspace{2mm}
\noindent
(i) $(H, [l,2l]_m)$.

\vspace{1mm}
\noindent
(ii) $(K,\{l\}_m)\cup (K^c, \{2l\}_m)\cup (H,[l+1,2l-1]_m)$. 

\vspace{1mm}
\noindent
(iii) $(K,\{l\}_m)\cup (K,\{2l+1\}_m)\cup (K^c,\{2l\}_m) \cup (K^c, \{l+1\}_m)\cup (H,[l+2,2l-1]_m)$. 

\vspace{2mm}
\noindent
Every sum-free subset of $G$ of the largest possible cardinality is one of the above for some splitting $(H,f)$ of $G$ by $\zm$ and some subgroup $K$ of $H$.
\end{thm}

\vspace{2mm}
\noindent
Given a splitting $(H,f)$ of $G$ by $\zm$ and a subgroup $K$ of $H$, with $G$ as in Theorem \ref{clas}, we shall hereafter write 
$L(H,f, 1), L(H,f,K,2), L(H,f,K,3)$ for the sets described respectively by $(i)$, $(ii)$, $(iii)$ of Theorem~\ref{clas}. Also, we shall say that a subset $L$ of $G$ is {\em presented by the splitting} $(H,f)$ if it is equal to one of these sets.

\vspace{2mm}
\noindent
To the extent we are aware, Theorem~\ref{clas} was formerly known only for groups of the form $({\Z}/p{\Z})^r$, due to Rhemtulla and Street~\cite{Street}, and in a small number of additional cases. We supplement Theorem 1.1 with the following result on the structure of sum-free subsets of $G$ of cardinality ``close" 
to the largest possible in $G$.

%\begin{defi}\label{3type}
%Let $G$ be a finite abelian group of type III.
%Given a splitting $(H,f)$ of $G$ by $\zm$ and a subgroup $K$ of $H$,
%we shall write $L(H,f, 1), L(H,f,K,2), L(H,f,K,3)$ to denote the sets given by $(i)$, $(ii)$, $(iii)$ of Theorem~\ref{clas} respectively. If a subset $L$ of $G$ is equal to one of these sets we say that $L$ is {\em presented by the splitting} $(H,f)$.
%\end{defi}

% Moreover, when this theorem is read together with the results of Diananda and Yap \cite{DY}, one obtains a complete characterisation of sum-free subsets of the largest cardinality in all (the three) types of finite abelian groups.

%\vspace{2mm}
%\noindent
%Let $A$ be a set from (i), (ii), or (iii) of Theorem~\ref{clas}. The
%first claim of Theorem~\ref{clas} is that $A$ is a sum-free subset of
%the largest possible cardinality in $G.$ That $A$ is sum-free  follows on noting that, given elements $(x_1,x_2), (y_1,y_2),
%(z_1,z_2) \in G$ with $(x_1,x_2) + (y_1, y_2) = (z_1,z_2)$, we must have  
% $x_1+y_1=z_1$ and $x_2+y_2=z_2.$ It is also easy to verify that $\Card(A) =
%\frac{ln}{m}$ which from~\eqref{mug} is equal to the
%maximum possible cardinality of a sum-free subset of $G$.

\vspace{3mm}
\begin{thm}\label{large}
Let $G$ be a finite abelian group of type III of exponent m and
cardinality 
$n$. For any sum-free subset $A$ of $G$ with ${\rm Card}(A) =
c(G) -\epsilon n$, where 
$\epsilon < {\rm min}(\frac{1}{6m}, 10^{-23}),$ there exists
a sum-free subset $L$ of the largest possible cardinality $c(G)$ in
$G$ such that ${\rm Card}(A\setminus L) \leq 4\epsilon n$.
\end{thm}

\vspace{2mm}
\noindent
The conclusion of Theorem~\ref{large} is essentially
best possible. Indeed, there are examples of sum-free sets $A$
satisfying the hypotheses of Theorem~\ref{large} such that there does
not exist any sum-free set $L$ of the largest possible cardinality in
$G$ with ${\rm Card}(A \setminus L)\leq \epsilon n.$

\vspace{2mm}
\noindent
In the second part of this article, comprising Sections 6 and 7, we give our first application of Theorem~\ref{clas}. More precisely, let ${\rm Aut}(G)$ and   $\mathcal{L}(G)$ respectively be 
the group of automorphisms and the set of all sum-free
subsets of the largest possible cardinality of a given finite abelian group $G$. Also, for any two subsets $A_1$ and $A_2$ of $G$ let us write $A_1 \sim A_2$ by $\operatorname{Aut}(G)$ if there is an $f$ in $\operatorname{Aut}(G)$ such that $f(A_1) = A_2$. The natural action of $\operatorname{Aut}(G)$ on the set of subsets of $G$ restricts to an action  of $\operatorname{Aut}(G)$
on $\mathcal{L}(G)$. One may therefore ask for a description of the orbit space   $\mathcal{L}(G)/ \operatorname{Aut}(G)$ under this action of $\operatorname{Aut}(G)$. Theorem \ref{sorb-BSP} below, which we obtain as an application of Theorem ~\ref{clas}, tells us that this question is  equivalent to the {\it Birkhoff subgroup embedding problem} when $G$ of type III. This problem asks to determine necessary and sufficient conditions so that $H_1 \sim H_2$ by $\operatorname{Aut}(G)$ where $G$ is a finite abelian group and $H_1$, $H_2$ are subgroups. We refer to  C.M. Ringel's lectures ~\cite{Ringel} for an exposition on Birkhoff's problem.

%We write $A_1 \sim A_2$ for any $A_1, A_2 \in \mathcal{L}(G)$
%when $A_1$ and $A_2$ are in the same orbit under this action of ${\rm Aut}(G)$ . 

%\vspace{2mm}
%\noindent
%Let $(H_1,K_1)$ and $(H_2,K_2)$ be pairs of abelian groups with $K_i$ a subgroup of $h_i$ for $i = 1,2$. We write  $(H,K_1)\sim (H,K_2)$ when there is an isomorphism of groups $f : H_1 \mapsto H_2$ such that $f(K_1) = K_2$. 
%If $f': H_1 \to H_2$ is an isomorphism, then it is easy to verify that the pairs $(H_1, K_1)$ and $(H_2,K_2)$ are isomorphic if and only if the pairs $(H_2, f'(K_1))$ and $(H_2, K_2)$ are equivalent. Since the structure theorem of finite abelian groups gives us a parametrisation of finite abelian groups up to isomorphism, it follows that to study the Birkhoff subgroup embedding problem, we may restrict to the case when $H_1=H_2$ and study the following problem.

%\begin{Q}[Birkhoff subgroup embedding problem]\cite{Birk} Let $H$ be a finite abelian groups. Then determine the equivalence classes of pairs $(H,K)$ with $K$ a subgroup of $H.$ \label{BSP}

%\end{Q}
%\noindent
%The readers interested in knowing the literature on this problem may  refer the lectures by C.M. Ringel~\cite{Ringel} and  the various references therein.

%\vspace{2mm}
%\noindent
%Using Theorem~\ref{clas}, we obtain the following result which shows that the %problems~\ref{sorb} and~\ref{BSP} are equivalent.

\begin{thm}\label{sorb-BSP}
Let $G$ be a finite abelian group of type III and exponent $m$ and let $(H,f)$ be a splitting of $G$ by $\zm$. Then every orbit of $\mathcal{L}(G)$ under $\operatorname{Aut}(G)$ contains an element of $\mathcal{L}(G)$ that is presented by $(H,f)$.
%Then given any $A\in \mathcal{G}$, there exists $A' \in \mathcal{L}(G)$ such that $A'$ has a presentation with respect to $(H,f)$ and $A\sim A'$. 
Moreover for $i, j\in \{2,3\}$, we have

\begin{equation}
L(H,f,K_1, i) \sim L(H,f,K_2,j) \;\text{by $\operatorname{Aut}(G)$} \iff K_1 \sim K_2 \; \text{by $\operatorname{Aut}(H)$ and  $i=j$.}
\end{equation}
\end{thm}

\vspace{2mm}
\noindent
The conclusions of Lemma~\ref{1split}, Propositions~\ref{orbp} and~\ref{orbp7} give some more information than provided by Theorem~\ref{sorb-BSP}. These results also discuss the orbit of $L(H,f,1)$ under $\operatorname{Aut}(G)$. As an easy consequence of  Theorem~\ref{sorb-BSP} we obtain :

%\vspace{2mm}
%\noindent
%Let $H$ be a supplement of a copy of $\zm$ in $G$ and
%$\mathcal{R}(H)$ be the set of subgroups of $H.$ Given $K \in \mathcal{R}(H)$ and $f \in {\rm Aut}(H)$, we have $f(K) \in \mathcal{R}(H).$ This defines an action of ${\rm Aut}(H)$ on~$\mathcal{R}(H).$
%By means of Theorem~\ref{sorb-BSP} we deduce fairly easily the following result.

\vspace{3mm}
\begin{thm}\label{rorb} Let $G$ be a finite abelian group of type III  of exponent $m$ and cardinality $n$ and let $H$ be a supplement in $G$ of a subgroup isomorphic to $\zm$. If $\mathcal{R}(H)$ is the set of subgroups of $H$ and  if $\mathcal{R}(H)/\operatorname{Aut}(H)$ is the set of orbits of $\mathcal{R}(H)$ under the natural action of $\operatorname{Aut}(H)$ then  we have  

\vspace{-3mm}
\begin{equation}
\label{orb}
{\rm Card}( \mathcal{L}(G)/{\rm Aut}(G)) = 2{\rm Card}( \mathcal{R}(H)/{\rm Aut}(H)) + \delta(m),
\end{equation}

\noindent
where $\delta(m)$ is 0 when $m=7$ and $1$ otherwise.
\end{thm}
\vspace{2mm}
\noindent
Theorem \ref{rorb} allows us to affirm the following  
 generalisation of a result of Rhemtulla and Street \cite{Street} for the groups $({\Z}/p{\Z})^r$, with $p$ a prime number congruent to 1 mod 3.

\vspace{3mm}
\begin{thm} \label{orbi} For integers $m$ each of whose prime divisors $p$ is congruent to $1$ mod $3$, the number of orbits under the action of the group of automorphisms of $({\Z}/m{\Z})^{r+1}$ on the set of sum-free subsets of the largest cardinality in $({\Z}/m{\Z})^{r+1}$ is 

\vspace{-3mm}
\begin{equation}
2\prod_{p|m} \binom{v_p(m) + r}{r} \; + \; \delta(m) \; , 
\end{equation}

\noindent
where $v_{p}(m)$ is the exponent of the prime $p$ in the prime decomposition of $m$  and $\delta(m)$ is $0$ when $m$ is $7$ and is $1$ otherwise.
\end{thm}

%\end{rem}
\vspace{2mm}
\noindent
Every subset of a sum-free subset is sum-free. Thus if   
${\rm SF}(G)$ is the set of all sum-free subsets of $G$ then ${\rm Card}({\rm SF}(G)) \geq 2^{c(G)}$, where, as before, $c(G)$ is the largest possible cardinality of a sum-free
subset of $G$. For $G$ of even cardinality an asymptotic formula for ${\rm Card}({\rm SF}(G))$ was obtained by Lev, Luczak and Schoen in~\cite{Lev} and independently by Sapozhenko in~\cite{Sapo}. A method for counting
sum-free sets was developed by Green and Ruzsa in a series of
papers. In~\cite[Theorem 1.9]{GR} they obtained an asymptotic formula
for ${\rm Card}({\rm SF}(G))$ in the case when the cardinality of $G$ is
divisible by a small prime divisor $q$ of the form $3k+2.$ 
However,  obtaining an asymptotic formula for ${\rm Card}({\rm SF}(G))$ when $G$ is of type III appears to be a rather difficult problem, even in the special case $G = ({\Z}/7{\Z})^r$ and in fact what is known are only upper and lower bounds for ${\rm Card}({\rm SF}(G))$ in terms of $2^{c(G)}$ and ${\rm Card}(G)$. 

\vspace{2mm}
\noindent
In \cite{GR} Green and Ruzsa show that for $G = ({\Z}/7{\Z})^r$ we have ${\rm Card}({\rm SF}(G)) \geq 2^{c(G) + c (\ln n)^2}$, where $n = 7^{r}$. Their  argument in fact yields a similar lower bound for all type III abelian groups $G$. For an upper bound, however, we only have ${\rm Card}({\rm SF}(G)) \leq 2^{c(G) + \frac{cn}{(\ln n)^{1/27}}}$ due to Balasubramanian and Gyan Prakash \cite{BGacta}. This is a slightly improved form of the original bound of Green and Ruzsa \cite{GR}, who had $\frac{1}{45}$ in place of $\frac{1}{27}$.  These upper bounds, however, hold for any finite abelian group~$G$, not necessarily of type III.

\vspace{2mm}
\noindent
In the third and final part of this article, comprising Sections 8 and 9, we apply 
 the  Theorem \ref{clas} to obtain the following result. 

\vspace{3mm}
\begin{thm}\label{coun} When $G$ is a finite abelian group of type III whose cardinality is $n$ and exponent $m$, the number of sum-free subsets of $G$ does not exceed $2^{c(G) + c_m \,n^{2/3}(\log n)^{4/3}}$, where $c(G)$ is the largest possible cardinality of a sum-free subset of $G$ and $c_m$ depends only on $m$.
\end{thm}

\vspace{2mm}
\noindent
For a finite abelian group $G$ whose cardinality is divisible by a prime $p \equiv 2$ mod 3, Green and Ruzsa obtained an asymptotic formula for $\Card(SF(G))$ by showing that the number of sum-free subsets of $G$ that are not contained in a sum-free subset of cardinality $c(G)$ is  $o_p(2^{c(G)})$. Recently Alon, Balogh, Morris and Samotij refined this result  in Theorem~1.1 of \cite{Alon} by showing that for each $m\geq c(q)\sqrt{n\log n}$ the number of sum-free subsets of $G$ of cardinality $m$ that are not contained in a sum-free of cardinality $c(G)$ is $o_p(\binom{c(G)}{m})$. 

\vspace{2mm}
\noindent
When $G$ is of type III, we are able to obtain only a result much weaker than those quoted in the preceding paragraph. Our result, given by Proposition~\ref{mpub}, is deduced from Theorem \ref{clas} via a modification of the arguments of Green and Ruzsa. We then use 
Proposition \ref{mpub} to obtain an apparently novel relation between the number of subsets with prescribed doubling in $H$, a supplement in $G$ of a subgroup isomorphic to ${\Z}/m{\Z}$, and the number of sum-free subsets of $G$. More precisely, for any positive integers $k_1$,$k_2$ and $H$ a finite abelian group, let  $S(k_1,k_2,H)$ be the number of subsets $B$ of $H$ with ${\rm Card}(B)=k_1$ and ${\rm Card}(2B) = k_2$. Then on setting 

\vspace{-3mm}
\begin{equation}
\label{ah}
a(H) \; =\; \sum_{k_1,k_2 \geq 1} \frac{S(k_1,k_2,H)}{2^{k_2}} \; 
\end{equation}

\noindent
we have the following theorem.
 
\begin{thm} \label{sfrelation} When $G$ is a finite abelian group of type III whose cardinality is $n$ and exponent $m$ and $H$ is a supplement of ${\Z}/m{\Z}$ in $G$ we have 

\vspace{-3mm}
\begin{equation}
\label{sk}
a(H) 2^{c(G)} \leq {\rm Card}({\rm SF}(G)) \leq n^2(a(H))^2 2^{c(G)} + o_{m}(2^{c(G)}) \; .
\end{equation}
\end{thm}

\noindent 
Theorem \ref{coun} is proved by combining \eqref{sk} with a bound for $\Card(S(k_1,k_2,H))$ obtained in \cite{cliquegen} generalising a bound for this quantity given by  Green\cite{Gclique} for the case when $H$ is of the form $({\Z}/p{\Z})^r$.

\vspace{2mm}
\noindent
Readers familiar with the works of Green and Ruzsa \cite{GR}  will recognise that our methods follow those of this work closely. We conclude this introduction by acknowledging our debt to these authors.

\section{Preliminaries}

\vspace{2mm}
\noindent
\subsection {Notation in Abelian Groups} We generally use $G$ to denote a finite abelian group. We use $A$,$B$, $C$, \ldots to denote subsets of $G$ and use $H$,$H^{\prime}$, $K$, $K^{\prime}$ for its subgroups. The group law in $G$ will be written additively with 0 for the identity element of $G$. Further, $n$ will  denote the cardinality of $G$ and $m$ its exponent. When $A$ and $B$ are subsets of $G$, $A+B$ will denote the image of the map $(x,y) \rightarrow x+y$ from $A \times B$ to $G$.

\begin{comment}
\vspace{2mm}
\noindent
{\em 2.2 Subsets of the ${\Z}$. ---} We use the usual notation for intervals in ${\bf R}$ to mean the set of integers contained in these intervals. When $X$ is a subset of  ${\Z}$ and $d \geq 1$ is an integer, $X_{d}$ will denote the canonical image in ${\Z}/d{\Z}$  of the set $X$. For example, when $a$ and $b$ are real numbers with $a < b$ and $d$ an integer $\geq 1$, $[a,b)$ denotes the set of integers in the real interval $[a,b)$ and $[a,b)_d$ denotes the canonical image of this set in ${\Z}/d{\Z}$.
\end{comment}

\subsection {Elementary Properties of Sum-free Subsets} Every subset of a sum-free subset of $G$ is a sum-free subset of $G$. The inverse image of a sum-free subset of $G$ under a homomorphism from $G^{\prime}$ to $G$ is a sum-free subset of $G^{\prime}$.

\subsection{Density}  When $G$ is a finite abelian group and $A$ is a subset of $G$ we write $\mu_{G}(A)$ to denote $\rm{Card}(A)/\rm{Card}(G)$. We call $\mu_{G}(A)$ the {\em density} of $A$. When $f$ is a surjective homomorphism of groups from $G$ onto $G^{\prime}$, we have the relation $\mu_{G}(f^{-1}(A^{\prime}))= \mu_{G^{\prime}}(A^{\prime})$ for every subset $A^{\prime}$ of $G^{\prime}$. For each integer $d \geq 1$, we write $\mu_d$ to denote $\mu_{{\Z}/d{\Z}}([d/3, 2d/3)_d)$. It is easily verified that 

\vspace{-3mm}
\begin{equation}
\mu_{d} = \frac{\left[\frac{d-2}{3}\right] + 1 }{d} \; .
\end{equation}

\vspace{2mm}
\noindent
\subsection{Density of Sum-free Subsets of the Largest Cardinality}
When $G$ is a finite abelian group we write $c(G)$ to denote the
cardinality of any sum-free subset of largest cardinality in $G$ and
write $\mu(G)$ to denote $c(G)/{\rm Card}(G)$. If $m$ is the exponent
of $G$, then for each divisor $d$ of $m$, ${\Z}/d{\Z}$ is a
quotient of $G$. It follows from the above that for
every divisor $d$ of $m$, $G$ contains a sum-free subset $A_{d}$ for
which $\mu_{G}(A_{d}) = \mu_{d}$. Consequently, $\mu(G) \geq \sup_{d |
  m} \mu_{d}$. The formula \eqref{mug} tells us that this inequality is
in fact an equality.

\vspace{2mm}
\noindent
\subsection{Types of Abelian Groups} \label{2.5}  A finite abelian group $G$ is said to be of  type I if  there exists a prime divisor of $n$ which is congruent to 2 modulo 3. When $G$ is of type I and if $p$ is the least among the primes congruent to 2 modulo 3 dividing $n$, we say that $G$ is of type I($p$). We say $G$ is of type~II if $G$ is not of type I and if $n$ is divisible by 3. Finally,  $G$ is said to be of type III if it is neither of type  I nor of type II . Thus $G$ is type III if and only if {\em all} prime divisors of $n$ are congruent to 1 modulo 3. With this division into three types, \eqref{mug} gives the following explicit relations for $\mu(G)$. 

\vspace{-2mm}
\begin{equation}
\label{mufinal} 
\mu(G) = 
\begin{cases}
\frac{1}{3} + \frac{1}{3p} \text{ when $G$ is of type I($p$)}, \\
\frac{1}{3} \text{ when $G$ is of type II}, \\
\frac{1}{3} - \frac{1}{3m} \text{ when $G$ is of type III.}
\end{cases}
\end{equation}

\subsection{A Consequence of the pigeonhole principle} We shall require the following lemma, which is an easy consequence of the pigeonhole principle.

\begin{lem}\label{pigeon}
Let $G$ be a finite group and $H$ be a subgroup. For some $x, y\in
G$, let $A$ and $B$ be subsets of $G$ with $A \subset H+x$ and $B
\subset H+y.$ If
$\min(\Card(A), \Card(B)) > \frac{\Card(H)}{2},$ then
$A+B = H+x+y.$
\end{lem}

\vspace{2mm}
\noindent
\subsection{An Application of Kneser's Theorem} Let $G$ be a finite abelian group acting on the set of its subsets by translation. For any $A\subset G$, we write $S(A)$ to denote the stabiliser of $A$ in $G.$ When $B$ and $C$ are subsets of $G$ such that ${\rm Card}(B+C)$ does not exceed ${\rm Card}(B) + {\rm Card}(C) -1$, we have by Kneser's theorem that 

\vspace{-2mm}
\begin{equation}
\label{kneser}
{\rm Card}(B+C) = {\rm Card}(B+H) + {\rm Card}(C+H) - {\rm Card}(H) \;,
\end{equation}

\noindent
where $H = S(B+C).$ It is now easily deduced that we have the following lemma.

\begin{lem}\label{Kneser}
Let $G$ be a finite abelian group and  $B \subset G$ with ${\rm Card}(B+B) < \frac{3}{2}{\rm Card}(B).$ Then $B+B$ is equal to a coset of $S(B+B)$ and consequently $B \subset S(B+B) +x$ for some $x\in G.$
\end{lem}

\subsection{Schur Triples} When $G$ is a finite abelian group of cardinality $n$ and $B$ is a subset of~$G$, an element $(x,y,z)$ of $B \times B \times B$ such that $x+y =z$ is called a {\em Schur triple}. We say that $B$ is an {\em almost sum-free} subset of $G$ if the number of Schur triples in $B \times B \times B$ is $o(n^2).$ The following results on almost sum-free subsets, due to Green and Ruzsa, will be crucial to our proof of the upper bound for ${\rm SF}(G)$ given by Theorem~\ref{coun}. 

\vspace{3mm}
\begin{thm}\cite[Theorem 1.5]{GSzemeredi} \label{almost} Let $G$ be a finite abelian group. Then every almost sum-free
 subset of $G$ may be written as $A \cup B$, where $A$ is sum-free and $
{\rm Card}(B)$ is $o(n)$.
\end{thm}

\vspace{3mm}
\begin{thm}\cite[Proposition 2.1']{GR} When the cardinality $n$ of a finite abelian group $G$ is sufficiently large, there is a family ${\mathcal F}$ of subsets of $G$ satisfying the following conditions. 
\label{granular}
\vspace{1mm}
\noindent
(i) Every sum-free subset of $G$ is contained in some element of the family ${\mathcal F}$. 

\vspace{1mm}
\noindent
(ii) There are no more than $2^{n(\log n)^{-1/18}}$ subsets of $G$ in ${\mathcal F}$.

\vspace{1mm}
\noindent
(iii) Every element of ${\mathcal F}$ is almost sum-free. In fact, ${\mathcal F}$ can be chosen so that the number of Schur triples in any element  of ${\mathcal F}$ does not exceed $n^2/(\log n)^{1/10}$.

\end{thm}

\vspace{2mm}
\noindent
\subsection{Sets with small sumset} Given a subset $B$ of an
abelian group $H$ and positive integers $k_1,k_2,$ recall that we write
$$S(k_1,k_2,H) =  \Card\left(\{B \subset H: \Card(B)= k_1, \Card(B+B) = k_2\}\right)$$
and
$$a(H) = \sum_{k_1,k_2}\frac{S(k_1,k_2,H)}{2^{k_2}}.$$

\vspace{2mm}
\noindent
When $H$ is a vector space over a finite field $\zp$, Ben Green has obtained
 an upper bound for the cardinality of $S(k_1,k_2,H)$ in~\cite[Proposition
 26]{Gclique}. In~\cite{cliquegen}, using a modification of
the arguments of Green,
the second author obtained an upper bound for $ \Card(S(k_1,k_2,H))$ for an
 arbitrary finite abelian group $H$ and proved the following
 result. 
\begin{thm}\cite[Theorem 6]{cliquegen}\label{smallG} Let $H$ be a finite abelian group of cardinality
$n$. Then the cardinality of $S(k_1,k_2,H)$ is at most
$$ n^{\frac{4k_2\log_2 k_1}{k_1}}\min(k_1^{c\omega(n)(k_1k_2\log k_1)^{1/3}}\binom{k_2}{k_1-1}(k_1^3 +1), k_1^{4k_1}),$$ where $\omega(n)$
denotes the number of distinct prime divisors of  $n$ and
$c$ is a positive absolute constant. 
\end{thm}

\vspace{2mm}
\noindent
 We shall also require the
following result from~\cite{GR}.
\begin{lem}\cite[Lemma 7.3 (ii)]{GR}\label{i2i}
Let $G$ be a finite abelian group of type III and $f: G \to \zq$ be a homomorphism. Then for any sum-free subset $A$ of $G$ and $i \in \zq$, we have
 $ \Card(A_i) +
 \Card(A_{2i}) \leq \frac{n}{q},$ where $A_i = f^{-1}\{i\} \cap A.$
\end{lem}
\begin{proof}
Since $A$ is sum-free,
 the set $A_i + A_i$ is disjoint from the set $A_{2i}$. Hence we have
$$\Card(A_{2i}) \leq \frac{n}{q}- \Card(A_i+A_i) \leq \frac{n}{q} - \Card(A_i).$$ 

\noindent
from which the lemma follows.
\end{proof}

\subsection{ The number of subsgroups of a finite group} We shall make use of the lemma below, which is a part of the folklore on finite groups. We  provide a proof for the lack of a suitable reference. In this subsection alone $G$ and $H$ denote arbitrary finite groups, not necessarily abelian.

\begin{lem}\label{nsub}
The number of subgroups of a finite group $G$ does not exceed~$|G|^{\log_2|G|}.$
\end{lem}

\begin{proof} Every finite group $G$ contains a generating set of cardinality not exceeding $\log_{2} |G|$. Indeed, supposing, as we may, that $|G| > 1$, let $X = \{g_1,\cdots,g_m\}$ be a generating set of $G$ satisfying the conditions $(i)$ $g_1$ is distinct from the identity element and  $(ii)$ $g_i$ does not belong to the subgroup generated by $\{g_1,\cdots, g_{i-1}\}$ for $2 \leq i \leq m$. Then we have that $2^m$ products $g_1^{n_1}\cdots g_m^{n_m}$, with $n_i \in \{0, 1\}$, are distinct elements of $G$ and this gives 
   
\begin{equation}\label{ubgen}
2^{\Card(X)} \leq |G|,
\end{equation}

\vspace{0.2cm}
\noindent
from which our assertion follows. Applying this to any subgroup of a given finite group $G$ of cardinality $n >1$, we see that every subgroup of $G$ is generated by a subset of $G$ of cardinality at most $\log_2 n$. Therefore the number of subgroups of $G$ is no more than 

$$\sum_{r=1}^{\lfloor \log_2n\rfloor}\binom{n}{r} \leq \sum_{r=1}^{\lfloor \log_2n\rfloor}\frac{n^r}{r!} \leq \frac{n^{\log_2n} \lfloor\log_2n\rfloor}{\lfloor \log_2 n\rfloor !} \leq n^{\log_2 n},$$ 

\noindent
as required. The conclusion of the lemma evidently holds when $|G| =1$ as well.
\end{proof}

\section{A Sketch of Proof of Theorems~\ref{clas} and \ref{large}.}
For the convenience of the reader we summarise here the proofs of Theorems~\ref{clas} and~\ref{large}. The principle is to obtain these results with the aid of the following
proposition and Kneser's theorem.

\begin{prop}\label{mainprop}
Let $G$ be a finite abelian group of type III, cardinality $n$ and exponent
$m.$ Let $A$ be a sum-free subset of $G$ with  $\Card(A) = c(G)
-\epsilon n$ with $\epsilon \leq \min(10^{-23},\frac{1}{6m}).$ Then
there exists a surjective homomorphism $f': G \to \zm$ such that 
\begin{equation}\label{slength}
A \subset f'^{-1}[2k,4k+1]_m,
\end{equation}
where $k= \frac{m-1}{6}.$ In other words there 
is a splitting $(H,f'')$
of $G$
by $\zm$ such that $H = \ker(f')$ and
\begin{equation}
A \subset (H, [2k,4k+1]_m).
\end{equation}

\end{prop}

\vspace{3mm}
\noindent
To prove Proposition~\ref{mainprop} we shall use the following result of
Green and Ruzsa.
\begin{prop}\cite[Propostion 7.2]{GR}\label{notor}
Let $G$ be a finite abelian group of type III and of cardinality $n$. Then
given any sum-free subset
$A$ of $G$ with $\Card(A) = c(G) - \epsilon n$, with $\epsilon \leq
10^{-23}$, there exists a surjective homomorphism $f: G \to
\Z/q\Z$ with $q \ne 1$ such that
\begin{equation}\label{llength}
A \subset f^{-1}[k+1,5k]_q,
\end{equation}
where $k = \frac{q-1}{6}$.
\end{prop}

\vspace{2mm}
\noindent
We provide a brief description of the arguments used in deducing
Proposition~\ref{mainprop} from Proposition~\ref{notor}. Let $A$ be as in Proposition~\ref{mainprop} and $f: G \to \Z/q\Z$ be a
surjective homomorphism as given by Proposition~\ref{notor}. Using the
arguments from~\cite{GR} and the assumed lower bound for the cardinality of $A$ it is easily seen that we have $q=m.$  To verify that the
assertion~\eqref{llength} may be strengthened to \eqref{slength}, we define a set $C(A) \subset \zm$ as follows:
\begin{equation*}
C(A) = \{i \in \zm: \Card(A\cap f^{-1}\{i\})> \frac{n}{2m}\}.
\end{equation*}
From the pigeonhole principle we deduce that $C(A)$ is a sum-free
subset of $\zm.$ An application of Kneser's theorem then shows that $\Card(C(A)) =
2k$; that is, in fact $C(A)$ is a sum-free subset of the largest possible
cardinality in $\zm.$ The
 structure of $C(A)$ is obtained by proving Theorem~\ref{clas} in the
 case when $G$ is cyclic. Using this we complete the proof of Proposition~\ref{mainprop}. Now for the details !

\section{Proof of Proposition~\ref{mainprop}}
%In this section we prove Proposition~\ref{mainprop}.
Let $A$ be as given in Proposition~\ref{mainprop}. Then from Proposition~\ref{notor} there exists a positive integer $q\ne 1$ and a surjective homomorphism $f: G \to \zq$ 
such that 
$$A\subset f^{-1}([k+1,5k]_q),$$
 where $k= \frac{q-1}{6}$.
 We shall first
observe that $q$ is equal to the exponent $m$ of $G$. 

\vspace{2mm}
\noindent
For any $i \in \zq$, we
write $A_i$ to denote the set $A \cap f^{-1}\{i\}$ and $\alpha_i$ to
denote the number $\frac{ \Card(A_i)q}{n}$. We write $I$ to denote the set
$[k+1,5k]_q$ and $I_0$ to denote the set $[2k+1, 4k]_q.$  
The following lemma was used in~\cite{GR} and is easy to check.
\begin{lem}\label{2k}
The set $I$ may be  divided into $2k$ disjoint pairs of elements of the
form $(y,y/2)$, with $y \in I_0$.
\end{lem}

\vspace{2mm}
\noindent
Since $A \subset f^{-1}(I)$, using Lemma~\ref{2k} we have 
\begin{equation}\label{i2ieqn}
\mu_{G}(A) = \frac{1}{q}\sum_{i
  \in \zm} \alpha_i  = \frac{1}{q}\sum_{i \in I_0}\left(\alpha_i +
  \alpha_{i/2}\right).
\end{equation}
The arguments used
to prove the following lemma are identical to that in~\cite{GR} to
deduce~\eqref{mug} from Proposition~\ref{notor} for type III groups. 
\begin{lem}\label{q=m} With the notations as above $q =m$.
\end{lem} 
\begin{proof}
 Using~\eqref{i2ieqn} and Lemma~\ref{i2i}, we obtain that $\mu_{G}(A) \leq  \frac{2k}{q} = \frac{1}{3} - \frac{1}{3q}$.
 Suppose the lemma is not true.
In that case, since $G$ is a type III group and $q$ divides $m$, we have $q \leq
\frac{m}{7}$. Therefore  we obtain that $\mu_{G}(A) \leq \frac{1}{3}
 -\frac{7}{3m} \leq \mu(G) - \frac{2}{m}$. This is contrary to the assumed lower bound of
$\mu_{G}(A)$. Hence we have $q=m$.
\end{proof}

\subsection{Reduction to cyclic case}
Given a set $A \subset G$ as in Proposition~\ref{mainprop}, we define $C(A)$
to be the  subset of the {\em cyclic group} $\zm$
as follows:
\begin{equation}
C(A) = \{i \in \zm: \alpha_i > \frac{1}{2}\}.
\end{equation}
In this
subsection we shall prove Proposition~\ref{reduction} stated below, which
shows that $C(A)$ is a sum-free subset of the largest cardinality in $\zm$.

\vspace{2mm}
\noindent
\begin{lem}\label{lsum-free}
For any $i, j \in C(A)$, we have $\alpha_{i+j} = 0$; in particular $C(A)$ is a
sum-free subset of $\zm$.
\end{lem}
\begin{proof}
Given any $i, j \in C(A)$, using Lemma~\ref{pigeon} with $H = \ker(f)$ we
obtain that $A_i + A_j = f^{-1}\{i+j\}$. Since $A$ is sum-free, 
the sets
$A_i + A_j$ and $A_{i+j}$ are disjoint. Hence the lemma follows.
\end{proof}
\begin{lem}\label{i2i>}
For any $i_0 \in [2k+1, 4k]_m$ we have $\alpha_{i_0} + \alpha_{i_0/2} \geq 1-
\left(\mu(G) - \mu_{G}(A)\right)m$.
\end{lem}
\begin{proof}
 Since $A \subset
f^{-1}(I)$, using~\eqref{i2ieqn} we have
\[
 \mu_{G}(A) \leq \frac{1}{m}\sum_{i \in I_0, i \neq
  i_0} (\alpha_i + \alpha_{i/2}) + \frac{1}{m}(\alpha_{i_0} + \alpha_{i_0/2}).
\]
Now using Lemma~\ref{i2i} we have that the first term in the right hand side of the above
inequality is at most $\frac{2k-1}{m} = \mu(G) -\frac{1}{m}$. Thus the result
follows.
\end{proof}
\begin{lem}\label{1/2}
For any $i_0  \in [2k+1,4k]_m$, exactly
one element from the pair of elements $(i_0,\frac{i_0}{2})$  belongs to $C(A)$.
\end{lem}
\begin{proof}
From Lemma~\ref{lsum-free}, $C(A)$ is sum-free. Therefore it can contain at most one element from the pair $(i_0,2i_0).$ Suppose the lemma is not true for some $i_0.$
 For the brevity of notation, let $i = i_0/2.$ Using Lemma~\ref{i2i>}, we obtain that
\begin{equation}\label{aia2i>}
\Card(A_i) + \Card(A_{2i}) > \frac{5n}{6m}.
\end{equation}
Since neither $i$ nor $2i$ is in $C(A)$, we obtain that 
$$\min( \Card(A_i),\Card(A_{2i})) > \frac{n}{3m}.$$
 Using the fact that the sets $A_i + A_i$ and $A_{2i}$ are
  disjoint subsets of $f^{-1}\{2i\}$ and~\eqref{aia2i>}, we also obtain
  that  
$$\Card(A_i + A_i) \leq \frac{n}{m} - \Card(A_{2i}) <
  \frac{n}{m}- (\frac{5n}{6m}- \Card(A_i))  < \frac{3}{2} \Card(A_i).$$
  Therefore applying Lemma~\ref{Kneser} with $\bbb = \ccc = A_i$, we obtain
  that $A_i +A_i = S(A_i+A_i)+g$, where $S(A_i+A_i)$ is the  stabiliser of $A_i + A_i$ in $G$ and $g\in G.$ Therefore we obtain that
$$ \Card(S(A_i+A_i)) =  \Card(A_i+A_i) \geq  \Card(A_i) > \frac{n}{3m}.$$
  We also have that $S(A_i+A_i)$ is a subgroup of $G$ contained in $H$ with $H= \ker(f).$ Since $H$ is a type III group, any proper subgroup of $H$ will have cardinality at most $\frac{n}{7m}.$ Hence we obtain that $S(A_i+A_i) = H$ and $A_i+A_i = f^{-1}\{2i\}.$ This implies that $A_{2i}= \emptyset$, which is contrary to our earlier conclusion that $ \Card(A_{2i}) > \frac{n}{3m},$ thus proving the lemma.
\end{proof}
Combining Lemmas~\ref{lsum-free} and~\ref{1/2}, we obtain the following result.
\begin{prop}\label{reduction}
The set $C(A)$ is a sum-free subset of $\zm$ with $ \Card(C(A)) = 2k$;
in other words, $C(A)$ is a sum-free subset of $\zm$ of the largest
possible cardinality. Further for any
$i, j \in C(A)$, we have $\alpha_{i+j} = 0$.
\end{prop}
\begin{rem}
We note that Lemma~\ref{1/2} also follows by appealing
to~\cite[Lemma 7.3 (iii)]{GR}. We have preferred to give a self contained proof.
\end{rem}
\subsection{Classification of sum-free subsets of the largest cardinality in cyclic groups.}
In this section, we prove Theorem~\ref{clas} in case $G = \zm.$ In particular, this  gives a structure of $C(A)$, when $A$ is a subset of a general finite abelian group. 

\vspace{2mm}
\noindent
Let $E$ be a sum-free subset of the largest cardinality in $\zm.$ From
Proposition~\ref{notor} and Lemma~\ref{q=m}, it follows that $d.E
\subset I = [k+1,5k]_m$ for some $d \in (\zm)^*,$ where $k=
\frac{m-1}{6}.$ Replacing $E$ by $d.E$, we may assume that $E \subset I.$ From~\eqref{mug}, we also know that $ \Card(E)=2k.$

\vspace{2mm}
\noindent
We write $E^c$ to denote the complement of $E$ in $\zm.$
For any subset $B$ of $\zm$ and an element $x \in \zm$, 
we write $\widetilde{B}$
and $\wtx$ respectively to denote their images in $\Z$ under the
natural unfolding map from $\zm$ to the interval $[0,m-1]$ in 
$\Z$. Notice that for any sets $B,C \subset \zm$, we have $\widetilde{(B\cap C)}= \widetilde{B}\cap \widetilde{C}.$

\vspace{2mm}
\noindent
We write $I_{-1}, I_0, I_1$ to denote, respectively, the subsets $[k+1,2k]_m$, $[2k+1,4k]_m$ and 
$[4k+1,5k]_m$  of $\zm$. For any set $B \subset \zm$ and $j \in \{-1,0,1\}$, we write $B_j,$ $\widetilde{B_j}$ to denote the sets $B \cap I_j,$ $\widetilde{B}\cap \widetilde{I}_j$ respectively.

\vspace{2mm}
\noindent
The following lemma is an easy consequence of Lemma~\ref{2k}. Noticing that $C(E)=E$, it may also be deduced from Lemma~\ref{1/2}.
\begin{lem}\label{xcx/2}
For any $x \in [2k+1,4k]_m$, we have
$$\frac{x}{2} \in E \iff x \in E^c.$$
\end{lem}
\begin{lem}\label{parity}
Let $x, y \in (E^c)_0$. If $\wt{x}$ and $\wt{y}$ are of same parity
{\em (}respectively of different parity{\em )}, then the element $\frac{x+y}{2}$
{\em (}respectively $\frac{x-y}{2}${\em )} belongs to $(E^c)_0$. 
\end{lem}
\begin{proof}
From the previous lemma we have $\frac{x}{2},
\frac{y}{2} \in E$. Since $E$ is sum-free we have that both the elements
$\frac{x+y}{2}$ and $\frac{x-y}{2}$ belong to $E^c$. If $\wt{x}$ and $\wt{y}$ are of
same parity then $\frac{x+y}{2}$ belongs to $I_0$ and hence to
$(E^c)_0$. Similarly if $\wt{x}$ and $\wt{y}$ are of different parity then
$\frac{x-y}{2}$ belongs to $I_0$ and hence to $(E^c)_0$. Hence the lemma 
follows.
\end{proof}
\begin{lem}\label{ap}
 When $\Card((\wtec)_0) \geq 2$, then $(\wtec)_0$ is an arithmetic progression and the common
  difference $d$
  between any two consecutive integers in  it is an odd integer. 
\end{lem}
\begin{proof}
Let $(\wtec)_0 = \{\wt{x_1} < \wt{x_2} <\cdots < \wt{x_t}\}$. Using Lemma~\ref{parity} it follows that for any $i$, the elements $\wt{x_i}$ and $\wt{x_{i+1}}$ are of different parity and 
$$\wt{x_i} = \frac{\wt{x_{i-1}}+\wt{x_{i+1}}}{2} \;\; \forall \; 2\leq i \leq t-1,$$
from which the lemma follows.
\end{proof}
The following lemma is easy to verify. 
\begin{lem}\label{x0x/2}
Let $x \in I_0$. If $\wt{x}$ is even, then $\frac{x}{2} \in I_{-1}$. If $\wt{x}$ is odd, then $\frac{x}{2} \in I_1.$ 
\end{lem}
On combining lemmas~\ref{xcx/2}, \ref{ap} and~\ref{x0x/2}, we obtain 
\begin{lem} \label{rele01} Let $E \subset \zm$ be as above. We then have
\begin{enumerate}
\item $\wte_{-1} = \{\frac{\wt{x}}{2}:  \wt{x} \in (\wtec)_0 \text{ and } \wt{x} \text{ is even.}\}.$
\item $\wte_{1} = \{\frac{m+\wt{x}}{2}:   
\wt{x} \in (\wtec)_0 \text{ and } \wt{x} \text{ is odd.}\}  $.
\item We have $$\Card((\wtec)_0) = \Card(\wt{E}_{-1})+\Card(\wt{E}_1).$$
\end{enumerate}
\end{lem}
\noindent  Lemmas~\ref{rele01} and \ref{ap} then give the following result.
\begin{lem}\label{ap-1}
 Let $i
  \in \{-1, 1\}$ and $ \Card(\wt{E}_i)\geq
    2$.  Then
  $\Card\left((\wtec)_0\right) \geq 2$ and $\wt{E}_i$ is an arithmetic progression. Moreover, the common difference of the arithmetic progression  
    $\wt{E}_i$  is 
    the same as the common difference of 
  the arithmetic progression $(\wtec)_0$. 
\end{lem}

\begin{lem}\label{2k0} Suppose  $m \neq 19$; that is we have $k \neq 3$.
If  the integer $2k$ does not belong to
$\wt{E}_{-1}$, then $\wt{E}_{-1} = \emptyset$. Similarly if  the integer $4k+1$ does not belong to
$\wt{E}_{1}$, then $\wt{E}_{1} = \emptyset$.
\end{lem}
\begin{proof}
It is sufficient to prove the claim for  $\wt{E}_{-1}$, since then for
$\wt{E}_{1}$, the claim follows replacing $E$ by $-E$.
The assertion follows trivially in case $k=1$. So we may assume that $k \geq 2$.

\vspace{2mm}
\noindent
Suppose the integer
$2k$ does not belong to $\wt{E}_{-1}$.
Using 
Lemma~\ref{rele01}, it follows that $4k \in
\wt{E}_0$. Since $E$ is sum-free and $2k-1 \con 4k + 4k \m $, it follows
that $2k-1$ does not belong to  $\wt{E}_{-1}$. In case $k =2$, we have
$I_{-1} = \{2k-1, 2k\}_m$ and hence the lemma follows in this case. Therefore
we are left to prove the lemma in the case when $k \geq 3$. 

\vspace{2mm}
\noindent
We claim that the set
$\wt{E}_{-1}$ does not contain any odd integer. Suppose the claim is not true and
$2k - 2r -1$ is the largest odd integer belonging to  $\wt{E}_{-1}$. Since we
know that $2k-1$ can not belong to $\wt{E}$, we have $r\geq 1$.  We have $\{2k
-2i -1: 0 \leq i \leq r-1\} \cup \{2k\}$ is a subset of $(\wt{E^c})_{-1}$ and hence using
Lemma~\ref{rele01}, $\{4k -2(2i+1): 0\leq i \leq r-1\} \cup \{4k\} \subset
\wt{E}_0$. Hence it follows that $\{8k -2i: 0\leq i \leq 2r-1\} \subset
\wt{E}_0 + \wt{E}_0$. Since $m=6k+1$ and $8k -2i \con 2k-2i -1 \m$, it follows that $\{2k
-2i -1: 0 \leq i \leq 2r-1\} \subset \wtec$.  Since $r \geq 1$, this implies that $2k-2r-1$ can not
belong to $\wt{E}$. Therefore it follows that $\wt{E}_{-1}$ does not
contain any odd integer.

\vspace{2mm}
\noindent
Using Lemma~\ref{ap-1}, it follows that $\Card(\wt{E}_{-1})\leq 1.$
Thus either $\wt{E}_{-1} = \emptyset$ or
$\wt{E}_{-1} = \{2k -2t\}$ for some $t\geq 1.$ To prove the lemma, we need to rule out the second possibility.

\vspace{2mm}
\noindent
Suppose $\wt{E}_{-1} = \{2k -2t\}$ with $t \geq 1$.  It follows
  using the
  Lemma~\ref{rele01} that the only even integer in $(\wtec)_0$ is $4k-4t$. 
  Since $E$ is sum-free we also have that $4k -(2k-2t) = 2k +2t$ belongs to  $(\wtec)_0.$ Therefore $2k+2t=4k-4t$ and hence $t = \frac{k}{3}$. 
 %We may note that the conclusion obtained so far holds also when $k =3$ and this will be used in the proof of the next lemma.

\vspace{2mm}
\noindent
Since we have $k \ne 3$, it follows that $t\neq 1$.
Therefore, we have that $2k+2\ne 2k+2t.$ Hence $2k+2$ belongs to $\wt{E}_0.$ Therefore the even integer $(2k+2)+(2k-2t)=4k-2t+2$ belongs to $(\wtec)_0$ and hence is equal to $4k-4t.$ This implies that $2t+2=0$, which is not possible as $t$ is a positive integer. Hence it
  follows that  $\wt{E}_{-1}$ is an empty set. Hence the
  lemma follows.
\end{proof}
\begin{lem}\label{2k19} Suppose $m = 19$ and hence $k =3$. If the integer $2k$ does not belong to $\wt{E}_{-1}$ then either $\wt{E}_{-1} = \emptyset$ or we have

\begin{eqnarray}
E & = & \{4\}_m \cup \{7,9,10,12\}_m \cup \{15\}_m \; \; \text{ and hence } \label{2k19e}\\
7.E & = & \{2k\}_m \cup [2k+2,4k-1]_m \cup \{4k+1\}_m. \label{2k19e2}
\end{eqnarray}
Similarly if the integer $4k+1$ does not belong to $\wt{E}_1$, then either $\wt{E}_1 = \emptyset$ or $E$ is equal to the set as in~\eqref{2k19e}.
\end{lem}
\begin{proof}
It is sufficient to prove the claim for  $\wt{E}_{-1}$, since then for
$\wt{E}_{1}$, the claim follows replacing $E$ by  $-E$. Now when $k=3$, we have $\tilde{I}_{-1} = \{4,5,6\} = \{2k-2,2k-1, 2k\},$ $\tilde{I}_0 = \{7,8,9,10,11,12\}$ and $\tilde{I}_1 = \{13,14,15\}.$ Suppose the integer
$2k$ does not belong to $\wt{E}_{-1}$.
Then as argued in Lemma~\ref{2k0}, using 
Lemma~\ref{rele01}, it follows that $4k \in
\wt{E}_0$. Since $E$ is sum-free and $2k-1 \con 4k + 4k \m $, it follows
that $2k-1 =5$ does not belong to  $\wt{E}_{-1}$. Therefore using Lemma~\ref{rele01}, it follows that $4k-2 = 10 \in \wte_0.$ Therefore if $\wte_{-1}\neq \emptyset$, then $\wt{E}_{-1} = \{4\}$ and  $\{10, 12\} \subset \wt{E}_0.$ This implies that $10+4 = 14 \in (\wtec)_{1}.$ Using Lemma~\ref{rele01}, it follows that $9 \in \wt{E}_0.$ This implies that 
$9 +4 = 13 \in (\wtec)_{1}$ which implies, using Lemma~\ref{rele01}, that $7 \in \wt{E}_0.$ This implies that $7+4 = 11 \in (\wtec)_0$, which using Lemma~\ref{rele01} implies that $15 \in \wt{E}_1.$ Therefore we obtain that $\{4,7,9,10, 12, 15\}_m \subset E.$ We verify easily that $\{4,7,9,10, 12, 15\}_m$ is a sum-free subset of $\zm$. When $m=19$, the cardinality of the largest sum-free set in $\Z/m\Z$ is equal to $6$, it follows that if $E_{-1} \neq \emptyset,$  then $E= \{4,7,9,10,12, 15\}_m$, which is same as in the right hand side of~\eqref{2k19e}.  Hence $7.E = \{28,49,63, 70,84,105\}_m 
= \{9, 11, 6, 13, 8, 10\}_m$ which is the same set as in~\eqref{2k19e2}. Hence the result follows.
\end{proof}

\begin{lem}\label{ec3}
Let $E \subset \zm$ be as above. In case $\Card\left((\wtec)_0\right) \geq 3$,
then $E = I_{-1} \cup I_1 = [k+1,2k]_m \cup [4k+1,5k]_m$; that is, $2.E = [2k+1,4k]_m$. 
\end{lem}
\begin{proof}
The lemma is equivalent to showing that  $(\wtec)_0 = [2k+1,4k].$ We prove this by showing that (i) $\{2k+1,4k\} \subset (\wtec)_0$ and (ii) the common difference $d$ between any two consecutive integer in $(\wtec)_0$ is equal to 1.

\vspace{2mm}
\noindent
From Lemmas~\ref{ap} and~\ref{rele01}, it follows that $\wt{E}_{-1} \ne \emptyset$ as well as $\wt{E}_1\ne \emptyset.$ Furthermore replacing $E$ by $-E$, if necessary, we may assume that $ \Card(\wt{E}_{-1})\geq 2.$ From Lemmas~\ref{2k0} and~\ref{2k19}, we have that $\{2k,4k+1\} \subset \wt{E}.$ This also implies that $\{2k+1, 4k\}
 \subset (\wtec)_0$. 

\vspace{2mm}
\noindent
Moreover, from Lemmas~\ref{ap}~and~\ref{ap-1}, the sets
 $\wt{E}_{-1}$ and $(\wtec)_0$ are arithmetic progressions with the same common
 difference $d$ which is an odd integer. In case $d \ne 1$, then $d \geq 3$
 and $\{2k+2, 2k+3\} \subset \wt{E}$. Since $d \geq 3$, we have
 $ \Card(I_{-1}) = k \geq 3$.  Let $2k-t$ be the second largest integer belonging to
 $\wt{E}_{-1}$. Then $t \geq 3$ and 
$\{(2k+2)+2k-t,(2k+3)+2k-t\}\subset (\wtec)_0$. This implies that $d=1$, which is contrary to the assumption  that
 $d\ne 1$. Hence the lemma follows.

\end{proof}
\noindent When $G = \zm$ and $(H,f)$ is a splitting of $G$ by $\zm$, then
evidently $H = \{0\}.$ Therefore when $G$ is cyclic
Theorem~\ref{clas} states the following.
\begin{thm}\label{cyclic}
Let $G = \zm$ be a type III group. Let $E \subset \zm$ be a sum-free
set with
$\Card(E) = 2k$, where $k = \frac{m-1}{6}$. Then for some $d \in (\zm)^*$ we
have that $d.E$ is one of the following three sets.
\begin{enumerate}
\item $[2k+1,4k]_m$.
\item $\{2k, 4k+1\}_m \cup [2k+2, 4k-1]_m$.
\item $[2k, 4k-1]_m$.
\end{enumerate}
\end{thm}
\begin{proof}
From Proposition~\ref{notor} and Lemma~\ref{q=m} we know that there exists $d
\in (\zm)^*$ such that $d.E \subset [k+1,5k]_m$. Replacing $E$ by $d.E$ we
assume that $d=1$.
The proof is divided into four cases according to the cardinality of
$(\wtec)_0$. In case $\Card\left((\wtec)_0\right) \geq 3$, then from
Lemma~\ref{ec3}, the set $2.E$ is equal to the set as in $(i)$ of the theorem. In
case $ \Card\left((\wtec)_0\right) =2$, then using Lemma~\ref{ap} and
Lemma~\ref{rele01} we have that $\Card(\wt{E}_{-1}) = \Card(\wt{E}_1) =
1$. If $m \neq 19$, then using Lemma~\ref{2k0} we obtain that $\wt{E}_{-1} = \{2k\}$ and
$\wt{E}_1 = \{4k+1\}$. Thus it follows that $E$ is as in $(ii)$ of the theorem. If $m = 19$, then using Lemma~\ref{2k19}, it follows that either $E$ or $7.E$ is as in $(ii)$ of the theorem. In
case $ \Card\left((\wtec)_0\right) =1$, then replacing $E$ by $-E$, if
necessary,  and using Lemma~\ref{rele01} we have $ \Card(\wt{E}_{-1}) =1$ and
$ \Card(\wt{E}_1) = 0$. Then using Lemmas~\ref{2k0} and~\ref{2k19}, we have $\wt{E}_{-1} =
\{2k\}$ and $\wt{E}_1 = \emptyset$. Thus it follows that $E$ is as in $(iii)$
of the theorem. In case  $\Card\left((\wtec)_0\right) =0$ we have trivially
that $E$ is as in $(i)$ of the theorem.
\end{proof}
Now using Theorem~\ref{cyclic} and Proposition~\ref{reduction}, we prove
Proposition~\ref{mainprop}.
\begin{proof}[Proof of Proposition~\ref{mainprop}]
Let $A$ be a set as in the proposition and $C(A)$ be the subset of $\zm$ as
above. From  Propositions~\ref{reduction} and~\ref{cyclic}, we have that there exists $d\in
(\zm)^*$ such that $d.C(A)$ is one of the three sets as given in
Proposition~\ref{cyclic}. We then verify that for any $i \notin
d^{-1}.[2k,4k+1]_m$, there always exist $i_1, i_2 \in C(A)$ such that $i = i_1 \pm
i_2$. Thus using Proposition~\ref{reduction}, we have $\alpha_{i} = 0$ for any
$i \notin d^{-1}[2k,4k+1]_m$. In other words 
\begin{equation}
A \subset (df)^{-1}[2k,4k+1]_m.
\end{equation}
Therefore~\eqref{slength} holds with $f' = df.$ Let $x \in f'^{-1}\{1\}$ and $f":
\zm \to G$ be the injective homomorphism satisfying $f''(1) = x.$ Then with $H = \ker(f')$ we have that
$(H,f'')$ is a splitting of $G$ by $\zm$ and 
\begin{equation}
A \subset (H,[2k,4k+1]_m).
\end{equation}
\end{proof}

\section{Proofs of Theorems~\ref{clas} and~\ref{large}}
Let $A$ be a set from (i), (ii), or (iii) of Theorem~\ref{clas}. The
first claim of Theorem~\ref{clas} is that $A$ is a sum-free subset of
the largest possible cardinality in $G.$ That $A$ is sum-free  follows on noting that, given elements $(x_1,x_2), (y_1,y_2),
(z_1,z_2) \in G$ with $(x_1,x_2) + (y_1, y_2) = (z_1,z_2)$, we must have  
 $x_1+y_1=z_1$ and $x_2+y_2=z_2.$ It is also easy to verify that $\Card(A) =
\frac{ln}{m}$ which from~\eqref{mug} is equal to the
maximum possible cardinality of a sum-free subset of $G$.

\vspace{0.5cm}
\noindent
In the rest of this section we prove the following result from which the second claim of Theorem~\ref{clas} as well as Theorem~\ref{large} are easily deduced. We recall that $k = \frac{l}{2} = \frac{m-1}{6}.$
\begin{prop}\label{xxl}
Let $A$ be as in Theorem~\ref{large}. Then there exists a splitting
$(H,f)$ of $G$ by $\zm$ and a subgroup $K$ of $H$, such that the following holds. With $L$ being
one of the following sets, $L(H,f,1), L(H,f,K,2), L(H,f,K,3)$,  we have
$ \Card(A\setminus L) \leq 4\epsilon n.$
\end{prop}

\vspace{2mm}
\noindent
Let $A$ be as in Proposition~\ref{xxl}. From Proposition~\ref{mainprop} there exists a splitting $(H,f')$ of $G$ by $\zm$ such that
$$A \subset (H, [2k,4k+1]_m).$$
It is easy to verify that $A \cup (H,[2k+1,4k-1]_m)$ is a sum-free subset of $G$. Without any loss of generality we may assume that $A$ is a maximal (with respect to set inclusion) sum-free set. Therefore $A$ is equal to the union of $(H,[2k+1,4k-1]_m)$ and
the following set:

\begin{equation}\label{2k,4k+1}
(A_{2k},\{2k\}_m)\cup (A_{2k+1},\{2k+1\}_m)\cup (A_{4k},\{4k\}_m)\cup (A_{4k+1},\{4k+1\}_m),
\end{equation}
with $A_{2k}, A_{2k+1}, A_{4k}, A_{4k+1}$ being subsets of $H.$
For any $i \in \{2k,4k+1\}_m$, applying Lemma~\ref{i2i>} with $i_0 = 2i$, we obtain the following inequality:

\begin{equation}
\Card(A_{i}) + \Card(A_{2i})  \geq  \frac{n}{m}- \epsilon n.\label{2k>}
\end{equation}

\begin{lem}\label{3/2}
If for some $i \in \{2k,4k+1\}_m$, we have $ \Card(A_i) > 2\epsilon n,$ then 
$ \Card(A_i+A_i) < \frac{3}{2}\Card(A_i).$
\end{lem}
\begin{proof}
Since $A$ is sum-free, we have $A_{2i}\subset H \setminus (A_i+A_i).$
Using this and~\eqref{2k>}, the lemma follows after a small calculation.
\end{proof}

Using Lemmas~\ref{3/2} and~\ref{Kneser}, we obtain the following result.
\begin{cor}\label{cstab}
If for some $i \in \{2k,4k+1\}_m$, we have $ \Card(A_i) > 2\epsilon n,$ then $A_i$ is contained in a coset of the stabiliser $K_i$ in $H$ of the subset $A_i+A_i$ of $H$; in other words, there exists an element $x_i \in H$ such that $A_i \subset K_i +x_i$ and $A_i+A_i = K_i +2x_i.$ 
\end{cor}

\begin{lem}\label{2k=-4k+1}
If $\min\left(\Card(A_{2k}), \Card(A_{4k+1})\right) > 2\epsilon n$, then 
$$A_{4k+1}+A_{4k+1} = -(A_{2k}+A_{2k}).$$
Moreover there exists a subgroup $K$ of $H$ and an element $x\in H$ such that
$A_{2k} \subset K +x$, $A_{4k+1}\subset K -x$, $A_{2k}+A_{2k}= K +2x$ and $A_{4k+1}+A_{4k+1}=K-2x.$
\end{lem}
\begin{proof}
Let $K_{2k}, K_{4k+1}$ be the subgroups and $x_{2k},x_{4k+1}$ be the elements in $H$ as given by Corollary~\ref{cstab}. To prove the lemma, we shall show that $K_{2k}=K_{4k+1}$ and $x_{2k}+x_{4k+1} \in K_{2k}.$ The lemma follows from this with the choice of $K=K_{2k}$ and $x=x_{2k}.$ 

\vspace{2mm}
\noindent
First  we prove the following facts: 
\begin{eqnarray}
A_{2k}-A_{4k} &=& H \setminus (K_{2k}-x_{2k})\label{2k-4k}\\
A_{4k+1}-A_{2k+1}   &=& H \setminus (K_{4k+1}-x_{4k+1}).\label{4k+1-2k+1}
\end{eqnarray}
Using the fact that $A_{4k} \subset H\setminus (K_{2k}+2x_{2k})$ and $A_{2k}\subset K_{2k} +x_{2k}$, it follows that 
\begin{equation}\label{sub}
A_{2k}-A_{4k} \subset H\setminus (K_{2k}-x_{2k}).
\end{equation}

\vspace{2mm}
\noindent
Since $A_{2k}\subset K_{2k}+x_{2k}$, we obtain that $ \Card(K_{2k}) > 2\epsilon n.$  Therefore it follows that $ \Card(A_{2k}) > \frac{\Card(K_{2k})}{2}.$ Let $y \in H$ be such that $K_{2k}+y \ne K_{2k}+2x_{2k}.$ Using~\eqref{2k>}, it follows that 
$ \Card(A_{4k} \cap (K_{2k}+y)) > \frac{\Card(K_{2k})}{2}.$ 

\vspace{2mm}
\noindent
Therefore using Lemma~\ref{pigeon}, it follows $(K_{2k}+x_{2k})-(K_{2k}+y)\subset A_{2k}-A_{4k}$ for any $y$ as above. Using this and~\eqref{sub}, we obtain~\eqref{2k-4k}. Using similar arguments, we obtain~\eqref{4k+1-2k+1}.

\vspace{2mm}
\noindent
Since $A$ is sum-free, it follows that $A_{2k}-A_{4k}$ and $A_{4k+1}$ are disjoint subsets of $H$. Hence we obtain from~\eqref{2k-4k} that $A_{4k+1} \subset K_{2k}-x_{2k}.$ Therefore 
$$K_{4k+1}+2x_{4k+1} = A_{4k+1}+A_{4k+1}\subset K_{2k}-2x_{2k}.$$
It follows that $K_{4k+1}\subset K_{2k}$ and $x_{2k}+x_{4k+1} \in K_{2k}.$ Similar arguments imply that $K_{2k}\subset K_{4k+1}.$ Hence the lemma follows.

\end{proof}

\vspace{2mm}
\noindent
The following lemma is easy to verify.
\begin{lem}\label{adsplit}
Let $(H,f')$ be a splitting of $G$ by $\zm$ and $x \in H.$ Further let $f: \zm \to G$ be the injective homomorphism with $f(2k) = f'(2k)+x.$ Then $H$ is a supplement of the image of $f$ in $G.$ Moreover given any $B \subset H$ and $\lambda \in \Z$ we have 
$$(B+\lambda x,\{\lambda 2k\}_m)_{(H,f')} = (B, \{\lambda 2k\}_m)_{(H,f)}.$$ 
\end{lem}
\begin{proof}[Proof of Proposition~\ref{xxl}] From Proposition~\ref{mainprop}, there exists a splitting $(H,f')$ of $G$ by $\zm$ such that $A$ is equal to the set as in~\eqref{2k,4k+1}. The proof of the proposition is divided into the following four cases.

\vspace{2mm}
\noindent
When ${\Card}(A_{2k}) \leq 2\epsilon n$ and $\Card(A_{4k+1})\leq 2\epsilon n:$ In this case with $L = L(H,f',0)$, we have 

$$\Card(A\setminus L) =  \Card(A_{2k}) + \Card(A_{4k+1})\leq 4\epsilon n.$$
Hence the proposition follows in this case.

\vspace{2mm}
\noindent
When $ \Card(A_{2k})> 2\epsilon n$ and $ \Card(A_{4k+1})\leq 2\epsilon n:$ Let $K_{2k}$ be a subgroup and $x_{2k} \in H$ be as in Lemma~\ref{cstab}. Let $f: \zm \to G$ be the injective homomorphism with $f(2k) = f'(2k) + x_{2k}.$ With $L = L(H,f,K,2)$, where $K= K_{2k}$, using Lemma~\ref{adsplit}, it follows that
$$\Card(A\setminus L) =  \Card(A_{4k+1}) \leq 2\epsilon n.$$
Hence the proposition follows in this case.

\vspace{2mm}
\noindent
When $ \Card(A_{2k})\leq 2\epsilon n$ and $ \Card(A_{4k+1})> 2\epsilon n:$ Replacing $(H,f')$ by $(H,-f')$, the proposition follows using the arguments of the previous case.

\vspace{2mm}
\noindent
When $\Card(A_{2k}) > 2\epsilon n$ and $ \Card(A_{4k+1}) > 2\epsilon n:$ Let $K$ be a subgroup and $x \in H$ be as in Lemma~\ref{2k=-4k+1}. Let $f: \zm \to G$ be the injective homomorphism with $f(2k) = f'(2k) +x.$ Then using Lemmas~\ref{2k=-4k+1} and~\ref{adsplit}, it follows that $A \subset L(H,f,K,3).$ Hence the proposition follows in this case.

\end{proof}

\section{Proof of Theorem~\ref{rorb}}\label{!present}
\noindent
Recall that $\llL(G)$ denotes the family of sum-free subsets of the largest
cardinality in $G$. We choose a splitting $(H,f)$ of $G$ by $\zm$ and write $\mathcal{R}(H)$ to denote the collection of subgroups of $H.$ We use $L,L_1,L_2,...$ to denote elements in $\llL(G)$ and use $K_1,K_2,...$ to denote elements in $\mathcal{R}(H).$ We recall that $k = \frac{l}{2} = \frac{m-1}{6}.$

\vspace{2mm}
\noindent
Given  $h \in {\rm Aut}(G)$ and for any $L \in \llL(G)$, we have that $h(L)$ also belongs to $\llL(G).$ This defines an action of ${\rm Aut}(G)$ on $\llL(G).$ Given $L_1,L_2 \in \llL(G)$, we say that $L_1 \sim L_2$  if $L_1$ and $L_2$ are in the same orbit.

\vspace{2mm}
\noindent
Given $h \in {\rm Aut}(H)$ and for any $K \in \mathcal{R}(H)$, we have that $h(K) \in \mathcal{R}(H).$ This defines an action of ${\rm Aut}(H)$ on $\mathcal{R}(H).$ We say $K_1 \sim K_2$, where $K_1,K_2 \in \mathcal{R}(H)$, if $K_1$ and $K_2$ are in the same orbit. In this section, we shall prove Theorem~\ref{rorb} which relates the number of orbits in $\llL(G)$ to the number of orbits in $\mathcal{R}(H).$

\vspace{2mm}
\noindent
We have the following two maps
\begin{equation}
T_1,T_2: \mathcal{R}(H)\to \llL(G)
\end{equation}
with $ T_1(K) = L(H,f,K,2)$ and 
$T_2(K) = L(H,f,K,3)$ for any $K \in \mathcal{R}(H)$.
We say that $L \in \llL(G)$ has a {\em presentation} with respect to $(H,f)$ if
$$L \in {\rm Im}(T_1) \cup {\rm Im}(T_2) \cup \{L(H,f,0)\}.$$
\begin{lem}\label{1split}
Given $L \in \llL(G),$ there exists $L_1 \in {\rm Im}(T_1) \cup {\rm Im}(T_2) \cup \{L(H,f,0)\} $ such that $L\sim L_1.$ 
\end{lem}
\begin{proof}
 From Theorem~\ref{clas}, there exists a splitting $(H',f')$ of $G$ by $\zm$ such that $L$ is one of the following sets $L(H',f',0), L(H',f',K',2), L(H',f',K',3)$ with $K'$ being a subgroup of $H'.$
It is easy to verify that there exists $h \in {\rm Aut}(G)$ such that $h(H') =H$, $h({\rm Im}(f')) = {\rm Im}(f)$ and $h(f'(1)) = f(1).$
 Then $h(L)$ is equal to one of the following sets: $L(H,f,0),\, T_1(h(K')), \,T_2(h(K')).$ Hence the lemma follows. 
\end{proof}
\noindent
In this section we shall prove the following results.
\begin{prop}\label{orbp} Let $G$ be a finite abelian group of type III and with exponent $m \ne 7$. Let $i \in \{1,2\}$ and $K_1,K_2 \in \mathcal{R}(H).$ Then the following hold:

\begin{eqnarray}
 T_i(K_1) &\sim& T_i(K_2) \iff K_1 \sim K_2,\\
T_1(K_1) &\nsim& T_2(K_2),\\
L(H,f,0) &\nsim& T_i(K_1).
\end{eqnarray}
\end{prop}

\begin{prop}\label{orbp7}
Let $G$ be a finite abelian group of type III and with exponent $m=7$; in other words, let $G = (\Z/7\Z)^r$. Let $i \in \{1,2\}$ and $K_1,K_2 \in \mathcal{R}(H).$ Then the following hold:

\begin{eqnarray}
T_i(K_1) &\sim & T_i(K_2) \iff K_1 \sim K_2,\\
T_1(K_1) &\nsim & T_2(K_2),\\
L(H,f,0) &\sim & T_i(K_1) \iff K_1 = H {\text and }\, i =2.
\end{eqnarray}

\end{prop}
\noindent Using Lemma~\ref{1split}, Propositions~\ref{orbp} and~\ref{orbp7}, it may be verified easily that Theorems~\ref{sorb-BSP} and~\ref{rorb} follows.

\begin{lem}\label{1imp}
Let $K_1,K_2 \in \mathcal{R}(H).$ Then for any $i \in \{1,2\},$ we have
 
$$K_1 \sim K_2 \implies T_i(K_1) \sim T_i(K_2).$$
 
\end{lem}
\begin{proof}
Since $K_1 \sim K_2$, there exists $h\in {\rm Aut}(H)$ with $h(K_1) = h(K_2).$
Moreover since $G = H\oplus {\rm Im}(f)$, we may extend $h$ to $\tilde{h} \in {\rm Aut}(G)$ by defining $\tilde{h}$ to be the identity map on ${\rm Im}(f).$ Then it is easy to verify that $\tilde{h}(T_i(K_1)) = T_i(K_2).$ Hence the lemma follows.
\end{proof}
Using Lemmas~\ref{1split} and~\ref{1imp} we obtain that
\begin{equation}
{\rm Card}( \mathcal{L}(G)/{\rm Aut}(G)) \leq 2{\rm Card}( \mathcal{R}(H)/{\rm Aut}(H)) + 1.
\end{equation}

\begin{comment}

\vspace{3mm}

\noindent
Recall that we have chosen a splitting $(H,f)$ of $G$ by $\zm$. This defines an isomorphism $\tilde{f}:H\oplus\zm\to~G$ defined as $\tilde{f}(x_1,x_2) = x_1 + f(x_2).$ For the purpose of the rest of this section, we shall identify $G$ with $H \oplus \zm$ with the aid of this particular isomorphism. 
 Given any
element $x \in G$, there exists a unique element $x_1 \in H$ and $x_2
\in \zm$ such that $x = x_1 + f(x_2):=(x_1,x_2).$ We have a natural projection map $\pi_2:H\oplus\zm \to \zm$ defined as $\pi_2((x_1,x_2)) = x_2$. As we have identified $G$ and $H\oplus \zm$, this defines a map
$\tilde{\pi}_2: G \to \zm$  with $\tilde{\pi}_2 = \pi_2\tilde{f}^{-1}(x)$, which, by an abuse of notation, we  continue to denote by $\pi_2$. Similarly, given an automorphism $h: G \to G$, the homomorphism $\pi_2h$ as well as $h\pi_2$ make sense as we have identified $\zm$ with its image in $G$ under the map~$f.$

\end{comment}

\vspace{2mm}
\noindent
Recall that for any subset $A$ of a finite abelian group $G$, the {\em stabiliser} $S(A)$ of
$A$ in $G$ is the subset of $G$ consisting of those elements $g\in G$ such
that $g+A=A$. Given any
element $x \in G$, there exists a unique element $x_1 \in H$ and $x_2
\in \zm$ such that $x = x_1 + f(x_2):=(x_1,x_2).$ Let $\pi_2: G \to \zm$ be the map given by $\pi_2(x) = \pi_2((x_1,x_2)) = x_2.$

\begin{lem}\label{unique}
Let $L$ be a set in $\llL(G)$ which has a presentation  with respect to $(H,f).$
 Then for any subset $A$ of $L$,  the stabiliser  $S(A)$ of $A$ in $G$ is
contained in $H$.
\end{lem}
\begin{proof}
The lemma is equivalent to showing that $\pi_2(S(A))=\{0\}$, which is equivalent to showing that $ \Card(\pi_2(S(A)))= 1.$
In case $ \Card(\pi_2(S(A))) = d$, then $\pi_2(S(A))$ consists of the image in
$\zm$ of all the integers  divisible by $\frac{m}{d}.$ Now clearly we have
that the set $\pi_2(A)$ is invariant under the translation by the elements in
$\pi_2(S(A))$. 
 Therefore if $t$ is any integer such that its residue modulo
$m$ belong to $\pi_2(A)$, then for any integer $i$, the residue modulo $m$ of
the integer $t + i\frac{m}{d}$
belongs to $\pi_2(A)$. Since $A\subset L$, we have $\pi_2(A) \subset
[2k,4k+1]_m$. Let $t$ be the largest integer in $[2k,4k+1]$ such that its
residue modulo $m$ belongs to $\pi_2(A)$. Since $t+
\frac{m}{d}$ also belongs to $\pi_2(A)$, it follows that
$$ t+ \frac{m}{d} \geq 2k+m = 8k+1.$$
In case $d \ne 1$, then since $G$ is of type III, we have $d$ is at least $7$.
Then since $t$ is at most $4k+1$, the left hand side of the above inequality is
at most $4k+1 + \frac{6k+1}{7}$, which is strictly less than the right hand
side of  the above inequality, which is absurd. Thus
$ \Card(\pi_2(S(A))=d=1$. Hence the lemma follows. 
\end{proof}
Using Lemma~\ref{unique}, the following result is easily obtained.
\begin{cor}\label{stab}
The stabiliser of $L(H,f,0)$ in $G$ is $H$ and
the stabiliser of $T_i(K)$ in $G$ is $K$ for any $i \in \{1,2\}.$ 
\end{cor}
We say that a subset $A$ of $G$ is  {\em almost translation invariant}
if it is invariant under translation by more than $\frac{n}{m}$ elements of
$G$; that is, if $ \Card(S(A)) \geq \frac{n}{m}.$ Using Lemma~\ref{unique}, we also obtain the following result.
\begin{cor}\label{not7alm}
Let $L \in  \llL(G)$ has a presentation with respect to $(H,f).$
The stabiliser of any nonempty almost translation invariant subset of $L$ is
$H.$ 
\end{cor}
\subsection{Proof of Proposition~\ref{orbp}}
In this subsection we shall assume that the exponent $m$ of $G$ is not $7.$
\begin{lem}\label{almnot7}
Given any $L \in \llL(G),$ there exists a nonempty almost translation invariant subset of $L.$ 
\end{lem}
\begin{proof}
Since $m$ is not $7$, the set $[2k+2,4k-1]_m$ is a nonempty set, where $k= \frac{m-1}{6}.$ Using this we verify that the set $(H,[2k+2,4k-1]_m)$ is a non-empty almost translation invariant subset of $L.$
\end{proof}
\begin{cor}\label{inv}
Let $L_1,L_2 \in \llL(G)$ have presentation with respect to $(H,f)$. If  $h \in {\rm Aut}(G)$ with $h(L_1) = L_2,$ then $h(H) = H.$
\end{cor}
\begin{proof}
From Lemma~\ref{almnot7}, there exists a nonempty almost translation invariant subset $B$ of $L_1$. From Corollary~\ref{not7alm}, the stabiliser of $B$ in $G$ is $H$. Since $h(L_1) = L_2$, it follows that $h(B)$ is an almost translation invariant subset of $L_2$ and the stabiliser of $h(B)$ in $G$ is $h(H).$ Using Corollary~\ref{not7alm}, it follows that $h(H) = H.$ Hence the claim follows.
\end{proof}
\noindent
Let $C_0, C_1, C_2$ are subsets of $\zm$ as follows. 
\begin{equation}\label{c0c1c2}
C_0 = [2k+1,4k]_m, \quad C_1= [2k,4k-1]_m, \quad C_2 = \{2k\}_m \cup
[2k+2,4k-1]_m \cup \{4k+1\}_m.
\end{equation}
\begin{lem}\label{cyclicorbit}
When $m \ne 7$, the sets $C_0, C_1, C_2$ lie in different orbits under the
action of $Aut(\zm)$. 
\end{lem}
\begin{proof}
It is easy to verify that $C_0 = -C_0$ and $C_2=-C_2$, whereas $C_1 \ne -C_1$. Therefore $C_1$
could  neither lie in the same orbit as $C_0$ nor it could lie in the same
orbit as $C_2$.  Using the assumption that $m\ne 7$, we may verify
that the cardinality of $C_2 +C_2$ is equal to $4k+1,$ whereas the cardinality of $C_0 +C_0$ is equal to $4k-1$. Hence $C_0$ could
not be in the same orbit as $C_2$. Thus the lemma follows.
\end{proof}

\begin{lem}\label{t1t2}
For any $i \in \{1,2\}$, and $K_1, K_2 \in \mathcal{R}(H)$, we have
$$T_i(K_1)\sim T_i(K_2) \implies K_1 \sim K_2.$$

\end{lem}
\begin{proof}
Let $h \in {\rm Aut}(G)$ with $h(T_i(K_1)) = T_i(K_2).$ To prove the lemma, we need to show that there exists $h' \in {\rm Aut}(H)$ with $h'(K_1) = h'(K_2).$
From Corollary~\ref{inv} we have $h(H) = H.$ Therefore the restriction $h'$ of $h$ to $H$ is an automorphism of $H.$ From Corollary~\ref{stab}, the stabiliser of $T_i(K_1)$ in $G$ is $K_1.$ Therefore it follows that $h(K_1)$ which is same as $h'(K_1)$ is the stabiliser of $T_i(K_2).$ But from Corollary~\ref{stab}, the stabiliser of $T_i(K_2)$ in $G$ is $K_2.$ Therefore $h'(K_1) = K_2.$ Hence the lemma follows.
\end{proof}
Let $\tp: G \to {\rm Im}(f) \subset G$ be the map defined as $\tp = f\pi_2.$
The following lemma is easy to verify.
\begin{lem}\label{resim}
Let $h \in {\rm Aut}(G)$ with $h(H) = H.$ Then the restriction of $\tp h$ to 
${\rm Im}(f)$ is an automorphism of ${\rm Im}(f).$ Moreover for any $A \subset G$, we have $ \tp h\tp(A) = \tp h(A).$
\end{lem}

\begin{lem}\label{t1nt2}
For any $K_1, K_2 \in \mathcal{R}(H)$, we have
$$T_1(K_1) \nsim T_2(K_2).$$
\end{lem}
\begin{proof}
Suppose the lemma is not true and there exist $K_1, K_2 \in \mathcal{R}(H)$ such that $T_1(K_1) \sim T_2(K_2).$ We claim that this implies that $\tp T_1(K_1)$ and $\tp T_2(K_2)$ are in the same orbit under the action of ${\rm Aut}({\rm Im}(f)).$

\vspace{1mm}
\noindent
There exists an $h \in {\rm Aut}(G)$ with $h(T_1(K_1))= T_2(K_2).$ Therefore we have $\tp h(T_1(K_1)) = \tp(T_2(K_2)).$
From Corollary~\ref{inv}, we have $h(H) = H.$  Using Lemma~\ref{resim} it follows that $\tp h\tp(T_1(K_1)) = \tp h(T_1(K_1)) = \tp T_2(K_2).$ Therefore the restriction of $\tp h$ to ${\rm Im}(f)$ is an automorphism of ${\rm Im}(f)$ which transports $\tp T_1(K_1)$ to $\tp T_2(K_2)$. Hence the claim follows.

\vspace{1mm}
\noindent
Unless $K_1 = K_2=H$, $ \Card(\tp (T_1(K_1))) \ne \Card(\tp(T_2(K_2)).$ Therefore it follows that $K_1 = K_2 = H.$ Therefore $\tp(T_1(K_1)) =f(C_1)$ and $\tp(T_2(K_2)) = f(C_2)$, where $C_1,C_2~\subset~\zm$ are defined as above. The lemma follows using Lemma~\ref{cyclicorbit}.
\end{proof}
Using arguments similar to those used in the proof of Lemma~\ref{t1nt2}, we obtain the following lemma.
\begin{lem}\label{0notti}
For any $K \in \mathcal{R}(H)$ and $i \in \{1,2\},$ we have
$$L(H,f,0) \nsim T_i(K).$$
\end{lem}
Using Lemmas~\ref{1imp},~\ref{t1t2},~\ref{t1nt2} and~\ref{0notti}, we obtain Proposition~\ref{orbp}.
\subsection{Proof of Proposition~\ref{orbp7}}
In this subsection, we shall assume that the exponent $m$ of $G$ is $7$; that is $G = (\Z/7\Z)^r.$  
\begin{lem}\label{alm}
Let $L \in \llL(G)$ have a presentation with respect to $(H,f).$
When $L \ne T_2(K)$ for any proper subgroup $K$ of $H$, then also there exists a nonempty almost translation invariant subset of $L.$ 
When $L = T_2(K)$ with $K$ being a proper subgroup of $H$, there does not exist any nonempty almost translation invariant subset of $L.$ 
\end{lem}
\begin{proof}
When $L \in {\rm Im}(T_1) \cup \{L(H,f,0)\}$, then $(H,\{3\}_7)$ is a nonempty almost translation invariant subset of $L$. When $L= T_2(H)$, then $L$ is a nonempty almost translation invariant subset of $L.$ This proves the first claim. Using Lemma~\ref{unique}, the second claim follows easily.
\end{proof}

\begin{lem}\label{alm7}
When $L = T_2(K)$ for some $K \in \mathcal{R}(H)$ with $K \ne H$, then  there exists a nonempty almost translation invariant subset of $L+L$. Further the stabiliser of any nonempty almost translation invariant
subset of $L+L$ in $G$ is  $H$.
\end{lem}
\begin{proof}
Since $K$ is a proper subgroup, we have that $ \Card(K) \leq
\frac{ \Card(H)}{7}$ and hence $\Card(K^c) \geq
\frac{6 \Card(H)}{7} > \frac{ \Card(H)}{2}$. Therefore $K^c + K^c = H$. Using this we verify that 
$$L+L = (H, [-1,1]_7)\cup (K^c,2)\cup (K,3)\cup (K,4),(K^c,5).$$
Therefore $(H,[-1,1]_7)$ is an almost translation invariant subset of $L+L.$ This
proves the first claim.

\vspace{2mm}
\noindent
We shall prove the second part of the lemma by showing that for any nonempty subset $A$ of $L+L$, we have
that $\pi_2(S(A))~=~\{0\}$. If not then 
for some subset $A$ of $L+L$ we have $\pi_2(S(A)) = \Z/7\Z$. This implies that
there are elements $x, y \in H$ such that  $(x, 1) \in S(A)$ and $(y,0) \in A$. 

\vspace{2mm}
\noindent
Then it follows that all the elements $(y+4x,4),
(y+5x,5), (y+3x,3)$ belong to $L + L$. This implies that the elements $y+4x$ and $y+3x$
belong to $K$ whereas the element $y+5x$ belongs to $K^c$. Thus the element
$y+5x -(y+4x) = x$ belongs to $K^c - K = K^c$. On the other hand $y+4x -
(y+3x)$ belongs to $K-K = K$. Thus the element $x$ belongs to the set $K$ as
well as its complement, which is not possible. Hence $\pi_2(S(A))=
\{0\}$. Hence the lemma follows.
\end{proof}
\begin{rem}
From Lemmas~\ref{alm}, \ref{alm7} and Corollary~\ref{not7alm}, the following result follows. If $L \in \llL(G)$ has a presentation with respect to splittings $(H_1,f_1)$ and $(H_2,f_2)$ of $G$ by $\zm$, then $H_1 = H_2.$
\end{rem}

\begin{lem}\label{cyc7}
Let $C_0 = [3,4]_7$, $C_1 = [2,3]_7$ and $C_2 = \{2,5\}_7.$
We have $C_0 = 2.C_2$ and $L(H,f,0) = 2.T_2(H).$ The sets $C_1$ and $C_0$ lie in different orbits under the action of ${\rm Aut}(\Z/7\Z).$
\end{lem}
\begin{proof}
The first claim is easy to verify. The second claim follows from the observation that $C_0=-C_0$, whereas $C_1\ne -C_1.$ Hence the lemma follows. 
\end{proof}

Using Lemmas~\ref{alm},~\ref{alm7},~\ref{cyc7} and arguments similar to those used in the proof of Proposition~\ref{orbp}, we obtain Proposition~\ref{orbp7}.

\section{Proof of Theorem~\ref{orbi}}
\noindent
Let $G = (\Z/m\Z)^{r+1}$, with every divisor of $m$ being congruent to $1$ modulo $3$, then $G$ is a finite abelian group of type III with $m$ being its exponent.
Let $H$ be a supplement of a copy of $\Z/m\Z$ in $G$, then $H = (\Z/m\Z)^r.$ In this section, we shall compute $ \Card(\mathcal{R}(H)/Aut(H))$ and use Theorem~\ref{rorb} to prove Theorem~\ref{orbi}.

\vspace{2mm}
\noindent
If $m = \prod_{p |m} p^{v_p(m)}$ then $H = \bigoplus_{p|m}H_p$ with $H_p =
(\Z/p^{v_p(m)}\Z)^r$ and we have
\begin{equation}\label{hp}
 \Card(\mathcal{R}(H)/Aut(H)) = \prod_{p} \Card(\mathcal{R}(H_p)/Aut(H_p)),
\end{equation}
where $\mathcal{R}(H_p)$ denotes the family of subgroups of $H_p$. Therefore we may
assume that $H$ is a $p$-group.

\vspace{2mm}
\noindent
A finite abelian $p$-group is called homogeneous of height $t$ and rank $r$ if it is isomorphic to the direct sum of $r$ copies of the cyclic group ${\Z}/{p^{t}}{\Z}$, where $t \geq 0$ and $r\geq 1$ are integers. 
Thus $H_p$ is a homogeneous group of height $v_p(m)$ and rank $r$.

\vspace{2mm}
\noindent
We will show that $\text{Aut}(H_p)$ acts transitively on isomorphism classes of subgroups of $H_p$. Thus if $K_1$ and $K_2$ are isomorphic subgroups of $H_p$ then there exists an automorphism of $H_p$ that transports $K_1$ onto $K_2$.

\vspace{2mm}
\noindent
We shall consider $H_p$ endowed with its natural ${\Z}$-module structure. Let $F$ be the free ${\Z}$-module  of rank $r$ and $M$ be the submodule $p^{v_p(m)} F$ of $F$. Then $M$ is free ${\Z}$-module of rank $r$ and $H_p$ is isomorphic to $F/M$. Let $f$ be an isomorphism from $F/M$ onto $H_p$ and let $\phi$ denote $f.p$, where $p$ is the canonical projection from $F$ onto $F/M$. Thus $\phi$ is a surjective homomorphism of ${\Z}$-modules from $F$ onto $H_p$ with $\text{Ker}(\phi) = M$.   

\begin{prop}
 Let $F$ be the free ${\Z}$-module  of rank $r$ and $M$ be the submodule $p^{v_p(m)} F$ of $F.$
When $E$ is a submodule of $F$ containing $M$ there exists a ${\Z}$-basis $\{e_1, e_2, \ldots, e_r\}$ for $F$ and an increasing sequence of integers $0 \leq a_1 \leq a_2 \leq \ldots \leq a_r \leq v_p(m) $ such that 
\begin{enumerate}
\item $\{p^{a_1}e_1, p^{a_2}e_2, \ldots, p^{a_r} e_r\}$ is a ${\Z}$-basis for $E$.
\item $E/M$ is isomorphic to $\bigoplus_{1 \leq i \leq r} {\Z}/p^{v_p(m)-a_i}{\Z}$.
\end{enumerate}
\end{prop}
\begin{proof}
Since $E$ contains $M$, $E$ is a submodule of $F$ of rank $r$. From the theory of modules over principal ideal domains it follows that there is a ${\Z}$-basis $\{e_1, e_2, \ldots, e_r\}$ of $F$ and an increasing sequence of integers $0 \leq n_1 \leq n_2 \leq \ldots \leq n_r$ such that $\{n_1e_1, n_2e_2, \ldots, n_r e_r\}$ is a ${\Z}$-basis for $E$. 

\vspace{2mm}
\noindent
Since $p^{v_p(m)} e_1 + p^{v_p(m)} e_2 +\ldots + p^{v_p(m)} e_r$ is in $M$ and therefore in $E$, there exist integers $c_i$, $1\leq i \leq r$, such that 

\vspace{-4mm}
\begin{equation}
p^{v_p(m)} e_1 + p^{v_p(m)} e_2 +\ldots + p^{v_p(m)} e_r = c_1 n_1 \,e_1 + c_2 n_2 \,e_2 + \ldots + c_r n_r\, e_r \; .
\end{equation}
\noindent
Equating the coefficients of the $e_i$ we deduce that, for each $i$, $1\leq i
\leq r$,  $n_i$ divides $p^{v_p(m)}$. On setting $n_i = p^{a_i}$ for each $i$,
$1\leq i \leq r$, we see that $\{a_i\}_{1 \leq i \leq r}$ is an increasing
sequence of integers in the interval $[0,v_p(m)]$ satisfying (i). Also (ii) is verified by noting that $M$ is generated by $\{p^{v_p(m)} e_1, p^{v_p(m)} e_2, \ldots, p^{v_p(m)} e_r\}$ and passing to quotients.
\end{proof}
\begin{prop}\label{homogorb}
If $K_1$ and $K_2$ are isomorphic subgroups of a homogeneous $H_p$ then there exists an automorphism of $H_p$ that transports $K_1$ onto $K_2$.
\end{prop}
\begin{proof}
 $K_1$ and $K_2$ are ${\Z}$-submodules of $H_p$ viewed as a ${\Z}$-module. Let $E_1$ and $E_2$ be the inverse images of $K_1$ and $K_2$ under $\phi$. Then $E_1$ and $E_2$ are submodules of $F$ containing $M$. Proposition~1 shows that there are bases $\{e_1, e_2, \ldots, e_r\}$ and $\{f_1, f_2, \ldots, f_r\}$ of $F$, increasing sequences $\{a_i\}_{1 \leq i \leq r}$ and $\{b_i\}_{1 \leq i \leq r}$ of integers in the interval $[0,k]$  such that $\{p^{a_1}e_1, p^{a_2}e_2, \ldots, p^{a_r} e_r\}$ is a ${\Z}$-basis for $E_1$ and $\{p^{b_1}f_1, p^{b_2}f_2, \ldots, p^{b_r} f_r\}$ is a ${\Z}$-basis for $E_2$. Moreover
\begin{equation}
\label{2}
\bigoplus_{1 \leq i \leq r} {\Z}/p^{k-a_i}{\Z} \,\cong \, E_1/M \, \cong \,K_1 \, \cong K_2 \, \cong \, E_2/M \, \cong \, \bigoplus_{1 \leq i \leq r} {\Z}/p^{k-b_i}{\Z} \;,
\end{equation}
\noindent
from which we have that $a_i = b_i$, for each $i$, ${1 \leq i \leq r}$. Thus
if  $\theta$ is the automorphism of $F$ defined by $\theta(e_i) = f_i$, for
each ${1 \leq i \leq r}$, then $\theta$ transports $E_1$ onto $E_2$ and leaves $M$
stable. On passing to quotients $\theta$ thus defines an automorphism of $H_p$
that transports $K_1$ onto $K_2$. 
\end{proof}
\noindent
One may easily verify the following lemma.
\begin{lem}\label{nosubg}
The number of isomorphism classes of subgroups of $H_p$ is equal to $\binom{v_p(m)+r}{r}.$
\end{lem}
Combining~\eqref{hp}, Proposition~\ref{homogorb} and Lemma~\ref{nosubg} we obtain.
\begin{prop}\label{orbhomog}
When $H = (\zm)^{r}$ then we have 
$$ \Card(\mathcal{R}(H)/Aut(H)) = \prod_{p^{v_p(m)}||m} \binom{v_p(m)+r}{r}.$$
\end{prop}
Combining Proposition~\ref{orbhomog} and Theorem~\ref{rorb}, Theorem~\ref{orbi} follows.
\section{Proof of Theorem~\ref{sfrelation}.} 
The main result we  require to prove Theorem~\ref{sfrelation} is the following result. The arguments used in proving it are similar to those used by Green and Ruzsa in proving~\cite[Lemma~5.8]{GR}.
\begin{prop}\label{mpub}
Let $G$ be a finite abelian group of type III and exponent $m$. With $o_m(2^{c(G)})$ exceptions, all sum-free $A \subset G$ are described as follows. Choose a splitting $(H,f)$ of $G$ by $\zm$ and take $A$ to be a subset of $(H,[2k,4k+1]_m)$, where $k= \frac{m-1}{6}.$
\end{prop}
Presently, we state and prove a few results required to deduce Theorem~\ref{sfrelation} from Proposition~\ref{mpub}. 
%The following lemma is easy to verify using Lemma~\ref{nsub}.
\begin{lem}\label{ubsplit}
The number of splittings $(H,f)$ of $G$ by $\zm$ is equal to $n(G)\frac{n}{m},$ where $n(G)$ is equal to the number of elements in $G$ of order $m.$
\end{lem}
\begin{proof}
Let $f: \zm \to G$ be an injective homomorphism, ${\rm Im}(f) = M$ and $f' : M \to \zm$ be the unique map such that $f'f$ is the identity map on $\zm.$ 
Let
\begin{equation}\label{split-ext}
F: \{H \subset G: H \oplus \zm = G\} \to \{\tilde{f}: \tilde{f}: G \to \zm, \; \tilde{f}|_M = f'\},
\end{equation}
with $F(H) = \tilde{f}$ be the unique homomorphism defined as $\tilde{f}|_H = 0$ and $\tilde{f}|_M = f'.$ Then we verify that the map $F$ is bijective by verifying that $F^{-1}(\tilde{f}) = {\rm Ker}(\tilde{f}).$  Therefore the number of $H$ such that $(H,f)$ is a splitting of $G$ is equal to the number of homomorphism $\tilde{f}$ such that its restriction to $M$ is equal to $f'.$

\vspace{0.2cm}
\noindent
Let $G = M \oplus H_0$, with $H_0$ being the direct sum of cyclic groups generated by $e_1,\cdots,e_r \in H_0$ with the order of $e_i$ being $n_i.$ Then the extension $\tilde{f}$ of $f'$ is uniquely determined by its value on $e_i.$ We may choose $\tilde{f}(e_i)$ to be any element $x \in \zm$ with $n_ix = 0$. Since $n_i$ divides $m$, the number of such $x$ is equal to $n_i$. Hence the number of extensions $\tilde{f}$ of $f'$ is equal to $\prod_{i=1}^rn_i = \frac{n}{m}.$ Therefore given any injective homomorphism $f: \zm \to G$, the number of $H$ such that $(H,f)$ is a splitting of $G$ by $\zm $ is equal to $\frac{n}{m}.$ Since the number of injective homomorphism $f$ is same as the number of elements in $G$ of order $m$, the result follows.
\end{proof}
\begin{cor}\label{nL}
When $G$ is a finite abelian group of type III which is of cardinality $n$ and
exponent $m$,
the total number of  sum-free subsets of the largest cardinality in  $G$
is at most $3n^22^{\log_2^2 n}.$
\end{cor}
\begin{proof}
From Theorem~\ref{clas}, any $L \in \llL(G)$ has a presentation with respect to some splitting $(H,f)$ of $G$ by $\zm$. The number of $L$ which have a presentation with respect to a splitting $(H,f)$ is no more than thrice the number of subgroups of $H$. 
Now from Lemma~\ref{nsub}, the number of subgroups of $H$ is at most  $2^{\log_2^2n}.$ Therefore using Lemma~\ref{ubsplit}, the lemma follows. 
 
\end{proof}

\begin{lem}\label{ubsf}
Let $(H,f)$ be a splitting of $G$ by $\zm.$ The number of sum-free subsets $A$ of $G$ with
$A\subset (H,[2k,4k+1]_m)$  is at most $a(H)^22^{c(G)}.$
\end{lem}
\begin{proof}
If $A\subset (H,[2k,4k+1]_m)$, then $A = \bigcup_{i\in [2k,4k+1]_m}(A_i,i)$ with the $A_i's$ being subsets of~$H.$ Since $A$ is sum-free, we have $A_{4k}\subset H\setminus (A_{2k}+A_{2k})$ and $A_{2k+1} \subset H\setminus (A_{4k+1}+A_{4k+1}).$ 
Therefore given any $B_1 \in S(k_1,k_2,H)$ and $B_2 \in S(k_1',k_2',H)$,  the number of $A's$, with $A_{2k}=B_1$ and $A_{4k+1}=B_2$,
is at most $\frac{1}{2^{k_2+k_2'}}2^{c(G)}.$ Using this and the definition of $a(H)$, the lemma follows.
\end{proof}
\begin{lem}\label{lb}
Let $(H,f)$ be a splitting of $G$ by $\zm.$ The number of sum-free subsets $A$ of $G$ with
$A\subset (H,[2k,4k]_m)$  is equal to $a(H)2^{c(G)}.$
\end{lem}
\begin{proof}
It is easy to verify that any $A=\bigcup_{i\in [2k,4k]_m}(A_i,i)$ with $A_{4k}\subset H\setminus (A_{2k}+A_{2k})$ is a sum-free subset of $G.$ For any $B \in S(k_1,k_2,H),$ the number of such $A's$ with $A_{2k}=B$ is equal to $\frac{1}{2^{k_2}}2^{c(G)}.$ Moreover for different choices of $B$, the sets $A$ are different. Hence the lemma follows. 
\end{proof}
\noindent
Combining Proposition~\ref{mpub} and Lemmas~\ref{ubsplit},~\ref{ubsf},~\ref{lb}, we obtain Theorem~\ref{sfrelation}.

\vspace{2mm}
\noindent
In the rest of this section we give a proof of Proposition~\ref{mpub}.The following result is an immediate corollary of Lemma~\ref{granular}.
%We say a subset $B$ of $G$ is {\em almost contained} in a subset $C$ of $G$ if
%$ \Card(B\setminus C) = o_m(n).$ Notice that if $B$ is almost contained in $C$
%and $C$ is almost contained in $D$, then $B$ is almost contained in $D$.
\begin{lem}\label{large-granular}
Let $\fff$ be a family of subsets as provided by Lemma~\ref{granular}. With $o(2^{\mu(G)n})$ 
exceptions, the rest of sum-free subsets $A$ of $G$ are subsets of some $F \in \fff$ with $\Card(F) \geq \mu(G)n - \frac{n}{(\log n)^{1/17}}.$
\end{lem}
Using Lemma~\ref{almost} and Theorem~\ref{large}, we easily verify the following result.
\begin{lem}\label{asubset}
Let $F$ be as in Lemma~\ref{large-granular}. For some sum-free subset $L$ of the largest cardinality in $G$, the set $F$ is almost contained in $L$; that is
 $F = B \cup
C$ with $B \subset L$ and $C\subset G$ with $ \Card(C) = o_m(n).$
\end{lem}
%\begin{proof}
%Let $\fff$ be a family of subsets of $G$ as provided by Lemma~\ref{granular}.
%Then it is clear that except $o(2^{\mu(G)n})$ sets, the rest of sum-free
%subsets $A$ are subsets of some $F_i \in \fff$ with $\Card(F_i) \geq \mu(G)n -
%\frac{n}{(\ln n)^{1/17}}$. Using this, Lemma~\ref{almost} and Theorem~\ref{large}, the lemma follows.
%\end{proof}

\noindent
Given any subset $L$ of $G$ and an element $x \in G$, let $R(x,L)$ be a
collection of pairwise disjoint two element subsets $B$ of $L$ such that $x$ is either a sum or
difference of elements in $B$. Moreover we assume that among all possible such
collections,  $R(x,L)$ is of largest possible cardinality.
The following lemma is easy to verify.
\begin{lem}\label{i1i2mid}
For any $i \in \zm$ with $i \notin [2k,4k+1]_m$, there exist $i_1,i_2 \in [2k+1,4k]_m$, with $i_1\ne i_2$ such that either $i = i_1+i_2$ or $i= i_1-i_2.$ When $k \neq 1$, similar conclusion holds with $i_1, i_2 \in [2k,4k-1]_m$ as well as with $i_1, i_2 \in \{2k, 4k+1\}_m \cup [2k+2,4k-1]_m.$
\end{lem}
\subsection{When $m$ is not $7$.}
In this section we shall prove Proposition~\ref{mpub} in the case when the exponent $m$ of $G$ is not $7.$     

\begin{lem}\label{fnot7}
Let the exponent $m$ of $G$ be distinct from $7$ and $L \in \llL(G).$ Let $(H,f)$ be a splitting of $G$ by $\zm$ such that $L$ has a presentation with respect to $(H,f).$
Then for any element $x\in G$ with
$x\notin (H,[2k,4k+1]_m)$, the cardinality of $R(x,L)$ is at least $\frac{5n}{7m}.$
\end{lem}
\begin{proof}
If $x \notin (H,[2k,4k+1]_m)$, then $x=(x_1,x_2)$ with $x_2 \notin [2k,4k+1]_m.$
We verify that there exists a proper subgroup $K$ of $G$ such that the following holds. Either  $(H\setminus K, [2k~+~1,4k]_m) \subset L$ or $(H, [2k,4k-1]_m) \subset L$ or $(H, \{2k, 4k+1\}_m) \cup (H, [2k+2, 4k-1]_m)\subset L.$  
Since $K$ is a proper subgroup of  the type III group $H$, we have $\Card(K) \leq \frac{\Card(H)}{7}.$
Using Lemma~\ref{i1i2mid}, it follows that 
$$\Card(R(x,L))\geq  \Card(H) -2\Card(K) \geq \frac{5n}{7m}.$$
 Hence the lemma follows.
\end{proof}

\begin{lem}\label{covnot7} Let $\fff$ be a family of subsets of $G$ as provided by Lemma~\ref{granular} and $F \in \fff$ with $\Card(F) \geq \mu(G)n - \frac{n}{(\log n)^{1/17}}.$ Let $L \in \llL(G)$ be a set such that $\Card(F\setminus L) = o_m(n)$ as provided by Lemma~\ref{asubset} and  $(H,f)$ be a splitting of $G$ by $\zm$ such that $L$ has a presentation with respect to $(H,f).$ Then the number of sum-free subsets $A\subset F$ 
 such that $A\not\subset~(H,[2k,4k+1]_m)$ is $O(2^{c(G)-c(m)n}),$ where $c(m)>0$ is a constant depending only upon $m.$
\end{lem}
\begin{proof}
Given any $x \notin (H,[2k,4k+1]_m)$, first we obtain an upper bound for the number of $A's$ containing $x.$ Since $A$ is sum-free, it may contain at most one element from the pair of elements $(y_1,y_2) \in R(x,L).$ Therefore we obtain that when $A$ contains $x$, the number of possibilities for $A\cap L$ is $2^{c(G) -2 \Card(R(x,L))}3^{ \Card(R(x,L))}.$ From Lemma~\ref{asubset}, we have $F = B \cup C$ and  $A\setminus L \subset C$ with $\Card(C) = o_m(n).$ Hence the number of possible subsets $A$ containing $x$ is at most
$$2^{o_m(n)}2^{c(G)-(2-\log3) \Card(R(x,L))}.$$
Since the total number of choices for $x$ is at most $n$, using Lemma~\ref{fnot7} we obtain the lemma.
\end{proof}
\noindent
Combining Lemmas~\ref{granular},~\ref{large-granular}, and \ref{covnot7}, we obtain Proposition~\ref{mpub} in the case when the exponent $m$ of $G$ is not $7$.

\subsection{When $m$ is $7$}
Now we prove Proposition~\ref{mpub} in case when $m$ is $7$; that is when $G = (\Z/7\Z)^r$ with $r$ being a positive integer. In this subsection, we shall assume that $m$ is $7.$ The following lemma is easy to verify.
\begin{lem}\label{s7}
Let $L \in \llL(G).$ Then there exists a splitting $(H,f)$ of $G$ by $\Z/7\Z$ such that $L$ has a presentation with respect to $(H,f)$ and $L \ne L(H,f,H, 3).$
\end{lem}
\begin{lem}\label{f7}Let $L \in \llL(G).$ Let $(H,f)$ be a splitting of $G$ by $\Z/7\Z$ as in Lemma~\ref{s7}.  Then for any nonzero $x\in G$ with $x \notin (H,[2,5]_7)$, the cardinality of $R(x,L)$ is at~least~$\frac{n}{49}.$
\end{lem}

\begin{proof}
When $L \ne (H,[2,3]_7)$, then $(H\setminus K, [3,4]_7) \subset L$ for some proper subgroup $K$ of $H$. 
 In this case the lemma follows using  arguments similar to those used in the proof of Lemma~\ref{fnot7}.
Now we prove the lemma, when
$L = (H, [2,3]_7).$
Let $x$ be a non zero element in $G$ with $x=(x_1, i)$ with $i \notin [2,5]_7$;
that is $i \in [-1,1]_7.$ The proof is divided into two sub-cases, when $i = 0$
and when $i \ne 0$.
 When $i \ne 0$ we can write $i$ either as a sum or
difference of two distinct elements in $[2,3]_7$. In this case, the lemma follows
using  arguments similar to those used in proving Lemma~\ref{fnot7}. To prove the
lemma when $i = 0$,
we observe that since $x \ne 0$, it follows that $x_1 \ne 0$. 
Therefore there exists a subgroup $H'$ of $H$ such that $x_1 \notin H'$ and
$ \Card(H')$ is at least $\frac{n}{49}.$ Therefore the sets $H'$ and $H'+x_1$ are disjoint.
Given any $z \in (H'+x_1, \{2\}_7)$ there exists a unique $y \in (H',
\{2\}_7)$ such that $x= z-y.$ Hence, in this case also the lemma follows.
\end{proof}
Using Lemma~\ref{f7}, we obtain that Lemma~\ref{covnot7} is also true when $m$ is $7.$ Using this, Lemmas~\ref{asubset} and \ref{nL}, we obtain Proposition~\ref{mpub} in the case when the exponent $m$ of $G$ is~$7$.
\section{Proof of Theorem~\ref{coun}}
Using Theorem~\ref{smallG} and the simple fact that $\binom{k_2}{k_1} \leq 2^{k_2},$ we verify the following fact.
\begin{lem}\label{bah}
Let $H$ be a finite abelian group of cardinality $n$ and exponent $m$. Let $a(H)$ be as defined in~\eqref{ah}. Then we have 
\begin{equation}
a(H) \leq n^{c\omega(m)n^{2/3}\log^{1/3} n},
\end{equation}
where $c$ is an absolute constant.
\end{lem} 
In proving the above lemma, we also use the fact that $\omega(n)=\omega(m).$
Combining Theorem~\ref{sfrelation} and Lemma~\ref{bah}, we obtain Theorem~\ref{coun}.

\section{Concluding remarks}
\begin{comment}
The results of this paper make use of Proposition~\ref{notor}, which was
proved
by Ben Green and Imre Ruzsa in~\cite{GR}. When there is a prime divisor of the
exponent of $G$ which is less than $1000$, the arguments of Green and Ruzsa
depend on highly tedious calculations which 
are done with the aid of computer. These calculations are not particularly
tedious when the smallest prime divisor of $m$ is sufficiently large.
Ben~Green has remarked to us that it would be highly desirable to obtain a different proof of this particular result of theirs.

\vspace{2mm}
\noindent
\end{comment}
The upper bound for $a(H)$ given by Lemma~\ref{bah} does not appear to be best possible. Any improvement will improve the result of Theorem~\ref{coun}. For any finite abelian group $H$ of cardinality $n$, one can show that $a(H)$ is greater than or equal to  the number of subgroups of $H$. Therefore when $H = (\Z/7\Z)^r$, we obtain that
\begin{equation}\label{lbah}
a(H) \geq 2^{c\log^2n},
\end{equation}
where $c>0$ is an absolute constant and $n=7^r$ is the cardinality of $H.$ Using Theorem~\ref{smallG}, one may also notice that the main contribution in the right hand side of~\eqref{ah} comes from those terms with $k_2$ close to $2k_1.$   

\vspace{2mm}
\noindent
Let $t\geq 2$ be a positive integer. We say that $A\subset G$ is $t$-free, if there is no solution of the equation $x_1+x_2+\ldots+x_t -y =0$ with $x_i's$ and $y\in A.$ We say that $G$ is of type $(t+1,1)$ if all the divisors of the exponent $m$ of $G$ are congruent to $1$ modulo $t+1.$ We write $\nu_t(G)$ to denote the number $\frac{\left[\frac{m-2}{t+1}\right]+1}{m}n.$ When $G$ is of type $(t+1,1)$, a conjecture of Hamidoune and Plagne~\cite{Plagne1}  states that the maximum possible cardinality $\mu_t(G)$ of $t$-free set in $G$ is equal to $\nu_t(G)$. 
 For any positive integers $t,k_1$,$k_2$, and $H$ any finite abelian group, we write $S(t,k_1,k_2,H)$ to denote the number of subsets $B$ of $H$ with ${\rm Card}(B)=k_1$ and ${\rm Card}(tB) = k_2$ and  set
 
 \vspace{-3mm}
\begin{equation}
a(t,H) \; =\; \sum_{k_1,k_2 \geq 1} \frac{S(t,k_1,k_2,H)}{2^{k_2}} \; .
\end{equation}
Using arguments similar to those used in the proof of Lemma~\ref{lb}, it follows that the number of $t$-free subsets in any finite abelian group $G$ of type $(t+1,1)$ is at least 
$$a(t,H)2^{\nu_t(G)},$$
where $H$ is a supplement of a copy of $\zm$ in $G.$

\vspace{3mm}
\noindent
{\bf Acknowledgement :} We are extremely grateful to the anonymous referees for carefully reading the manuscript and providing us with a number of useful comments. We are thankful to Professor Amritanshu Prasad for drawing our attention to the Birkhoff subgroup embedding problem and related literature.
%\bibliographystyle{plain}
%\bibliography{master}

%\begin{comment}
\vspace{1cm}
\begin{flushleft}

{\em
Institute of Mathematical Sciences, \hfill Harish-Chandra Research Institute,\\
C.I.T. Campus, Taramani, \hfill Chhatnag Road, Jhunsi,\\
Chennai - 600113, India. \hfill Allahabad - 211 019, India.\\
e-mail: balu@imsc.res.in \hfill e-mail: gyan.jp@gmail.com, gyan@hri.res.in
  \\
\hfill suri@hri.res.in}

\end{flushleft}
%\end{comment}

\end{document}